\documentclass[a4paper,11pt, reqno]{amsart}
\usepackage{array}
\usepackage{graphicx}
\usepackage{float}
\usepackage{amsmath,latexsym}
\usepackage{amssymb}
\usepackage{geometry}
\usepackage{amsfonts}
\usepackage{hyperref}
\usepackage{color}
\usepackage{multirow}

\usepackage{natbib}

\usepackage{setspace}
\allowdisplaybreaks

\let\origmaketitle\maketitle
\def\maketitle{
  \begingroup
  \def\uppercasenonmath##1{} 
  \let\MakeUppercase\relax 
  \origmaketitle
  \endgroup
}

\usepackage{graphicx}%
\usepackage{multirow}%
\usepackage{amsmath,amssymb,amsfonts}%
\usepackage{amsthm}%
\usepackage{mathrsfs}%
\usepackage[title]{appendix}%
\usepackage{xcolor}%
\usepackage{textcomp}%
\usepackage{manyfoot}%
\usepackage{booktabs}%
\usepackage{algorithm}%
\usepackage{algorithmicx}%
\usepackage{algpseudocode}%
\usepackage{listings}%
\usepackage{color}%
\usepackage{rotating}  
\usepackage{float}
\usepackage{adjustbox}

\usepackage{xcolor}
\definecolor{armygreen}{rgb}{0.19, 0.53, 0.43}

\title[]{\large A Mathematical Programming Approach to Optimal Classification Forests }

\author[V. Blanco, A. Jap\'on, J. Puerto \MakeLowercase{and} P. Zhang]{{\large V\'ictor Blanco$^\dagger$, Alberto Jap\'on$^{\ddagger,\star}$, Justo Puerto$^\ddagger$, and Peter Zhang$^\star$}\medskip\\
$^\dagger$Institute of Mathematics (IMAG), Universidad de Granada\\
$^\ddagger$Institute of Mathematics (IMUS), Universidad de Sevilla\\
$^\star$ Heinz College of Information Systems and Public Policy, Carnegie Mellon University}

\date{\today}

\begin{document}

\maketitle

\begin{abstract}
This paper introduces Weighted Optimal Classification Forests (WOCFs), a new family of classifiers that takes advantage of an optimal ensemble of decision trees to derive accurate and interpretable classifiers. We propose a novel mathematical optimization-based methodology which simultaneously constructs a given number of trees, each of them providing a predicted class for the observations in the feature space. The classification rule is derived by assigning to each observation its most frequently predicted class among the trees. 
We provide a mixed integer linear programming formulation (MIP) for the problem and several novel MIP strengthening / scaling techniques. We report the results of our computational experiments, from which we conclude that our method has equal or superior performance compared with state-of-the-art tree-based classification methods for small to medium-sized instances. We also present three real-world case studies showing that our methodology has very interesting implications in terms of interpretability. Overall, WOCFs complement existing methods such as CART, Optimal Classification Trees, Random Forests and XGBoost. In addition to its Pareto improvement on accuracy and interpretability, we also see unique properties emerging in terms of different trees focusing on different feature variables. This provides nontrivial improvement in interpretability and usability of the trained model in terms of counterfactual explanation. Thus, despite the apparent computational challenge of WOCFs that limit the size of the problems that can be efficiently solved with current MIP, this is an important research direction that can lead to qualitatively different insights for researchers and complement the toolbox of practitioners for high stakes problems.

\end{abstract}

\section{Introduction}
When it comes to creating a machine learning model, two fundamental factors are considered. On the one hand, the aim is to create models capable of making accurate predictions. On the other hand, we seek interpretability, since interpretable classifiers are preferable to black-box models (see e.g., \cite{IML}). Both features are known to be compatible although there often exists a trade-off between them. Accurate models tend to be more complex, and therefore difficult to interpret, whereas interpretable ones are usually less accurate. This is the case,  for example, with Classification Trees (CTs) and Random Forests (RFs). \\

CTs are a family of classification methods built on a hierarchical relation among a set of nodes. A tree involves a set of branches connecting the nodes and defining the paths that observations can take by following a tree-shaped scheme. At the first stage of a CT method, all the observations belong to one node, which is known as the root node. From this node, branches are sequentially created by splits on the feature space, creating intermediate nodes (branch nodes) until the terminal nodes (leaf nodes) are reached. The predicted label for an observation is given by the most frequent class of the leaf node where it is located. \\
 
In contrast to the interpretable CT, \cite{RF} introduced RF a method that creates an ensemble of trees (also called a forest) in order to combine their predictions so as to make a final prediction by means of voting. In RF, the number of trees involved in the ensemble is fixed, and usually large, which makes the model less interpretable.  Each tree in the forest is constructed by drawing a random sample from the original training dataset and applying an algorithm to construct a classification tree. In general, RFs are more accurate models than a single decision tree, and more stable in terms of the selected points involved in the training dataset \citep{RF,stable_class,gov22}.  Furthermore, in case the tree ensemble has a small number of trees, it provides a flexible classification rule that allows a simpler analysis of counterfactual explanations by combining the interpretation of the paths given by the trees. Determining the classification of an observation within a particular class provided that the classifier assigned to it a different class can be done in many different ways by obtaining the majority of the votes among all the trees. However, finding the trade off between accuracy and the complexity of the ensemble (in terms of number of trees and their depths) is a challenging question that we address in this paper.\\

CTs were initially presented from an algorithmic point of view, normally taking into account greedy heuristic principles \citep{CART}. However, while Optimal Classification Trees (OCTs)
have been widely studied in the literature in recent years and have proven to provide numerous advantages (see e.g., \citep{OCT_survey}), to the best of our knowledge, the question of finding an \emph{optimal} forest has not been previously studied. \\

In this paper, we take a first step in this direction and study the (Weighted) Optimal Classification Forest (WOCF). We provide a mixed integer linear programming (MILP) formulation for the binary classification version of the problem. We show in our computational experiments that compared with a deep tree, a WOCF can improve accuracy in small to medium-sized datasets with just a small number of shallow trees. This indicates that WOCFs can potentially provide Pareto improvements over current methods, on both accuracy and interpretability for small to medium-sized instances. Moreover, we illustrate by means of some real examples how the solutions obtained are not only easier to interpret than those obtained with a CT, but also more flexible when being used. This is due to the fact that the first splits made in a deep CT have a strong impact on the following splits to be applied, whereas in an WOCF this effect is mitigated due to the independence of the splits distributed among the different trees.

\subsection{Literature}

There are many different heuristic approaches to building a CT. The most popular method is CART, introduced by \citet{CART}. CART is a greedy heuristic approach that myopically constructs a binary CT without anticipating future splits. Starting at the root node, it constructs the splits by means of separating hyperplanes minimizing an impurity measure at each of the branch nodes. Each split results in two new nodes and this procedure is repeated until a stopping criterion is reached. The trees derived from CART follow a top-down greedy approach, a strategy that is also shared by many other popular methods like C4.5~\citep{quinlan93} or ID3~\citep{quinlan96}. One of the advantages of these top-down approaches is that the trees can be obtained rapidly even for large datasets, as the whole process is based on solving manageable problems at each node. Nevertheless, the optimality of these trees is not guaranteed because at each node the best split is sought without taking into account the splits to be made in deeper nodes. According to this, the obtained tree structures may not capture the correct distribution of the data, which may result in out-of-sample errors when making predictions. Moreover, these methods can generate very deep (complex) trees if high training set accuracy is desired, and this results in a direct loss of interpretability as well as an increased risk of overfitting. This problem is usually overcome by applying a post-process pruning of the trees, which generally consists of analyzing if it is worth splitting a given node by its improvement in the local classification rate. Different approaches to construct \emph{optimal} learning trees have been already proposed, as the one in \cite{hu2020learning}, where the MAXSAT framework is applied.\\

On the other hand, \cite{OCT} recently proposed an optimization approach to build a CT. This approach has been shown to outperform heuristic approaches. Optimal Classification Trees (OCTs) are capable of capturing hidden splits that may not have a strong local impact on the training task but allow the final model to be more accurate. OCTs have been an active field of study in recent years due to their good performance, and several algorithms have been designed to efficiently train them. \cite{demirovic2022}, \cite{lin2020} and \cite{hu2019optimal} propose different dynamic programming based algorithms to construct OCTs. Constraint Programming and SAT have also contributed to the OCT field with works such as those of \cite{verhaeghe2020learning}, \cite{ yu2020computing}, \cite{hu2020learning} or \cite{narodytska2018learning}. Moreover, we can see in \citep{verwer2019learning} an alternative integer programming formulation to construct OCTs with a smaller number of decision variables than the one introduced by \cite{OCT}, as well as a column generation solution approach by \cite{firat2020column}.\\

In order to obtain easily interpretable CTs, these are usually derived in their univariate version, i.e., with splits in which only one predictive variable is involved. Nevertheless, oblique CTs, trees where multiple predictive variables can be involved in the splits, have also been widely studied and one can find numerous heuristic methods in the literature~\citep{oblique_ct,sparse_o_trees}, as well as the oblique version of OCT, (OCT-H), introduced in \cite{OCT}, or other optimal methods as the one presented in \cite{MOCTSVM}.\\ 

As for the field of RF, since its introduction by \cite{RF}, it has become one of the most popular Machine Learning ensemble models, and some others such as gradient-boosted decision tree methods have been inspired by it. RF is a method that combines bagging \citep{bagging} and CART strategies, where different training sub-samples are obtained to grow CTs so that one can combine their predictions by means of a voting strategy. Its popularity is due to the fact that it has great advantages in terms of accuracy and training stability with respect to CTs,  to its good performance for default hyperparameters, and because it can be easily adapted to parallel computational strategies due to its own nature. Therefore, improvements in this method with different ways of obtaining the training samples and ways to randomize the trees have been studied in the literature \citep[see e.g.][]{aglin2021,bernard2009,RFsurvey}. 
Furthermore, an implementation of boosted trees, derived from a trees ensemble, {\color{blue} namely XGBoost,} has become one of the most popular machine learning methods as we can see reflected, for instance, in the fact that it was used in most of the winning ideas in Kaggle data science competitions in 2015 \citep{xgboost}. However, despite its many advantages, the major drawback of the method is the loss of interpretability as the number of trees in the ensemble and the depth of trees increase.\\

In recent years, mathematical optimization has brought great advances in Machine Learning. As well as in tree-shaped classification rules, in methods such as Support Vector Machines (SVM), mathematical optimization plays a crucial role as the classifiers are derived directly from the solution of an optimization problem \citep{svm,Multiclass,lpSVM}. Furthermore, by means of properly modifying the objective function and constraints of a given optimization problem, one can adapt the methods to scenarios such as variable selection \citep{gunlunk,baldomero2020tightening,baldomero2021robust,Amorosi2024}, problems with special accuracy requirements \citep{benitez2019cost,gan2021robust}, or problems with unbalanced and noisy instances \citep{eitrich2006efficient,OCTSVM,blanquero2021optimal}. For this reason, numerous studies combining machine learning and mathematical optimization continue to be carried out and are a field of great interest because of their successes in practical applications.\\

Nevertheless, as far as we know, there is not much work in the literature on the use of mathematical optimization to derive tree ensembles. \cite{o_treeensembles} considers a given tree ensemble that predicts a target continuous variable using controllable independent variables (given a certain application), and it studies how to set these independent variables in order to maximize the predicted value by means of a MILP. Apart from this, we could not find in the literature anything regarding optimal tree ensembles for classification problems.

\subsection{Contributions and Structure}

This paper is among the first to study the concept of Weighting Optimal Classification Forest (WOCF). 

\begin{enumerate}
    \item We propose a novel MIP model that formulates the WOCF problem. It can naturally and modularly port different formulations of the recently proposed Optimal Classification Tree models, while adding some nontrivial generalizations that consider the combinatorial properties (enumeration, weighting, and voting) of multiple trees. 
    \item We design novel MIP strengthening cuts that break symmetry and improve the formulation since trees are permutation-invariant in this voting forest. We show in our numerical experiences and case studies that this WOCF method can do well on problems with $10^2$ to $10^4$ observations and can handle from 10 to 100 variables.
    \item For practitioners, WOCF complements the existing methods (CART, OCT, RF) and provides Pareto improvement on accuracy and interpretability for small to medium sized problems. This makes it especially suited for high stakes problems with moderate datasets. The main advantage of WOCFs as opposed to existing RF is to drastically reduce the number of trees involved in the ensemble, resulting in more interpretable classification rules without significant loss of accuracy. Across 20 test datasets, we show a fair comparison in terms of accuracy using WOCFs, CART, OCTs, RFs, and XGB showing that WOCF results match or outperform other existing CT-based methods and obtain comparable results with respect to RFs for small and medium-sized instances. Furthermore, we analyze in detail three real-world case studies where we further analyze the solutions obtained by WOCFs, concluding that our tool provides classifications rules which are easy to interpret and flexible to apply, as opposed to more rigid and complex solutions provided by complex CTs.
    \item We see unique properties emerging from the trained models of WOCF, compared with previous methods. In particular, we see each tree being delegated a different ``perspective'' in the classification and voting process, leading to each tree focusing on a different subset of variables and different subregions of the feature space. This is somewhat surprising to us at first, but not unreasonable.  Indeed, one could imagine an \emph{optimally} constructed forest would have a much better combinatorial structure than \emph{random} forests. At the same time, an optimally constructed \emph{forest} could have a much more efficient and cleaner way to use distinct stumps / shallow trees to cover the feature set and feature space than an optimally constructed \emph{tree} that has to rely on deeper and less efficient cuts with inter-dependence among splitting variables. For research purpose, we think this is an interesting future direction to study the combinatorial behavior of optimal forests. For practical purpose, this not only provides nontrivial improvement in interpretability, but also usability of the trained model in terms of counterfactual explanation.
\end{enumerate}

The idea of WOCF is not surprising given the recent development of Optimal Classification Trees over traditional CARTs. But we believe the results are somewhat unexpected and qualitatively different as mentioned in the previous list of contributions, and therefore we believe this line of research holds significant potential despite its apparent present day computational limitation for large instances (this could change quickly). It can already be a useful tool for practitioners for medium sized problems when accuracy and interpretability are important. It shows very interesting signs of new combinatorial properties for machine learning and optimization researchers.

In Section \ref{sec:prel} we introduce the problem under study and fix the notation used throughout the paper. Section \ref{sec:form} presents the mathematical optimization model for the binary WOCF. In Section \ref{sec:exps} we report the results of our computational experiments comparing WOCFs with RFs and CT-based methodologies. In Section \ref{sec:case} we analyze in detail the solutions provided by WOCFs and OCTs over three real case studies. Finally, in Section \ref{sec:conc} we report our conclusions.

\section{Preliminaries}\label{sec:prel}

In this section, we introduce the concept of WOCF and recall the main elements that set the base for the approach that is presented in Section \ref{sec:form}, i.e. the main concepts regarding CTs. 

All throughout this paper, we consider that we are given a training sample $\mathcal{X} =\{ (x_1,y_1), \ldots , (x_n,y_n) \}  \subseteq \mathbb{R}^p \times \left\lbrace 0, 1 \right\rbrace$, where $p$ features have been measured for a set of $n$ individuals ($x_1, \ldots, x_n$) as well as a label in $\{0, 1\}$ for each of them $(y_1,\ldots , y_n)$. The goal of binary supervised classification is to derive, from $\mathcal{X}$, a decision rule $D_\mathcal{X}: \mathbb{R}^p \rightarrow \{0, 1\}$ capable of accurately predicting the label of out-of-sample observations. We assume, without loss of generality, that the features are normalized, i.e., $x_1, \ldots, x_n \in [0,1]^p$, and we denote by $N=\{1, \ldots, n\}$ the index set for the observations in the training dataset.\\

As already pointed out, CTs are classifiers based on a hierarchical relationship among a set of nodes. The decision rule for a CT method is built by recursively partitioning the feature space by means of hyperplanes. In the first stage, a root node for the tree is considered which all the observations belong to. Branches are sequentially created, at the so-called branch nodes, by splits on the feature space, creating intermediate branch nodes until a leaf node, nodes where no more splits are made, is reached. Specifically, at each branch node, $t$, of the tree a hyperplane $\mathcal{H}_t = \{z \in \left[0,1\right]^p:  {(a_t)}^\intercal z = b_{t} \}$ is constructed (note that  ${(a_t)}^\intercal $ stands for the transpose operator applied to the vector $a_t\in \left[0,1\right]^p$)  and the splits are defined as $ {(a_t)}^\intercal z  <  b_{t} $ (left branch) and $ {(a_t)}^\intercal z \geq  b_{t} $ (right branch).  In a CT, we call depth the maximum number of splits that can be made before reaching a leaf node. Given a maximum depth, $D$,  a CT can have at most  $T=2^{D+1} -1$ nodes. These nodes are differentiated into two types:
\begin{itemize}
\item Branch nodes:  $\tau_B = \left\lbrace 1, \ldots , \lfloor T/2\rfloor \right\rbrace$ are the nodes where splits can be applied.
\item Leaf nodes: $\tau_L = \left\lbrace  \lceil T/2\rceil ,\ldots, T \right\rbrace $ are the nodes where predictions for observations are performed.
\end{itemize}

On the other hand, according to the notation used to describe the hierarchical relationships in the tree, we define $p(t)$ as the parent of node $t$, for $t=2, \ldots, T$. We define
 $A_{L}(t)$ (resp $A_{R}(t)$) as the set of left (resp. right) ancestors of node $t \in \{2, \ldots, T\}$, which are the nodes whose left (resp right) branch has been followed on the path from the root node to $t$. \\

\begin{figure}[h!]
    \centering
    \includegraphics[width=0.82\textwidth]{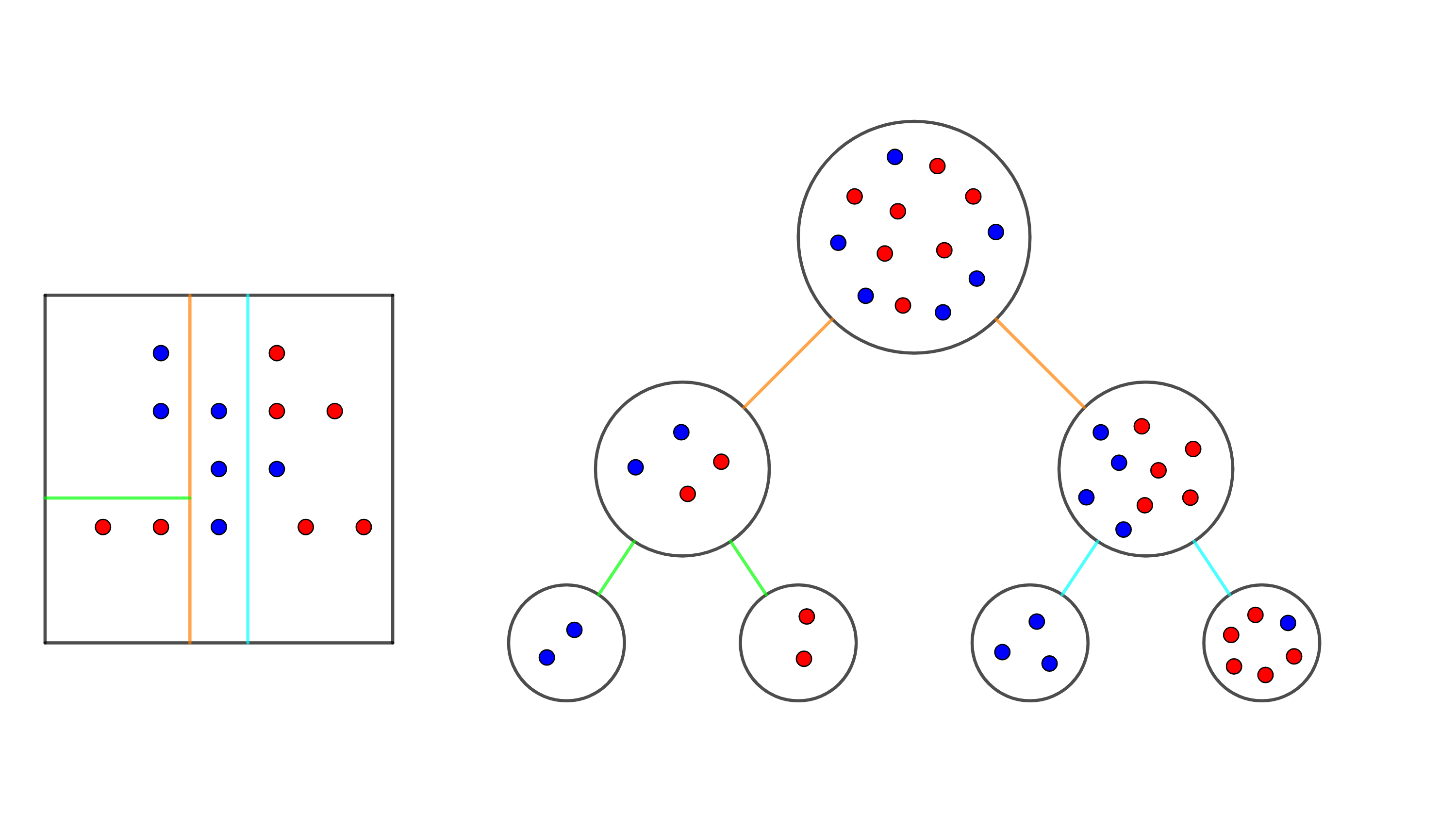}
    \caption{Binary classification problem with its split solution in $\left[0,1\right]\times [0,1]$ (left) and its OCT graph-solution form (right)}
    \label{fig:1}
\end{figure}

\begin{figure}[h!]
\centering
\includegraphics[width=0.5\textwidth]{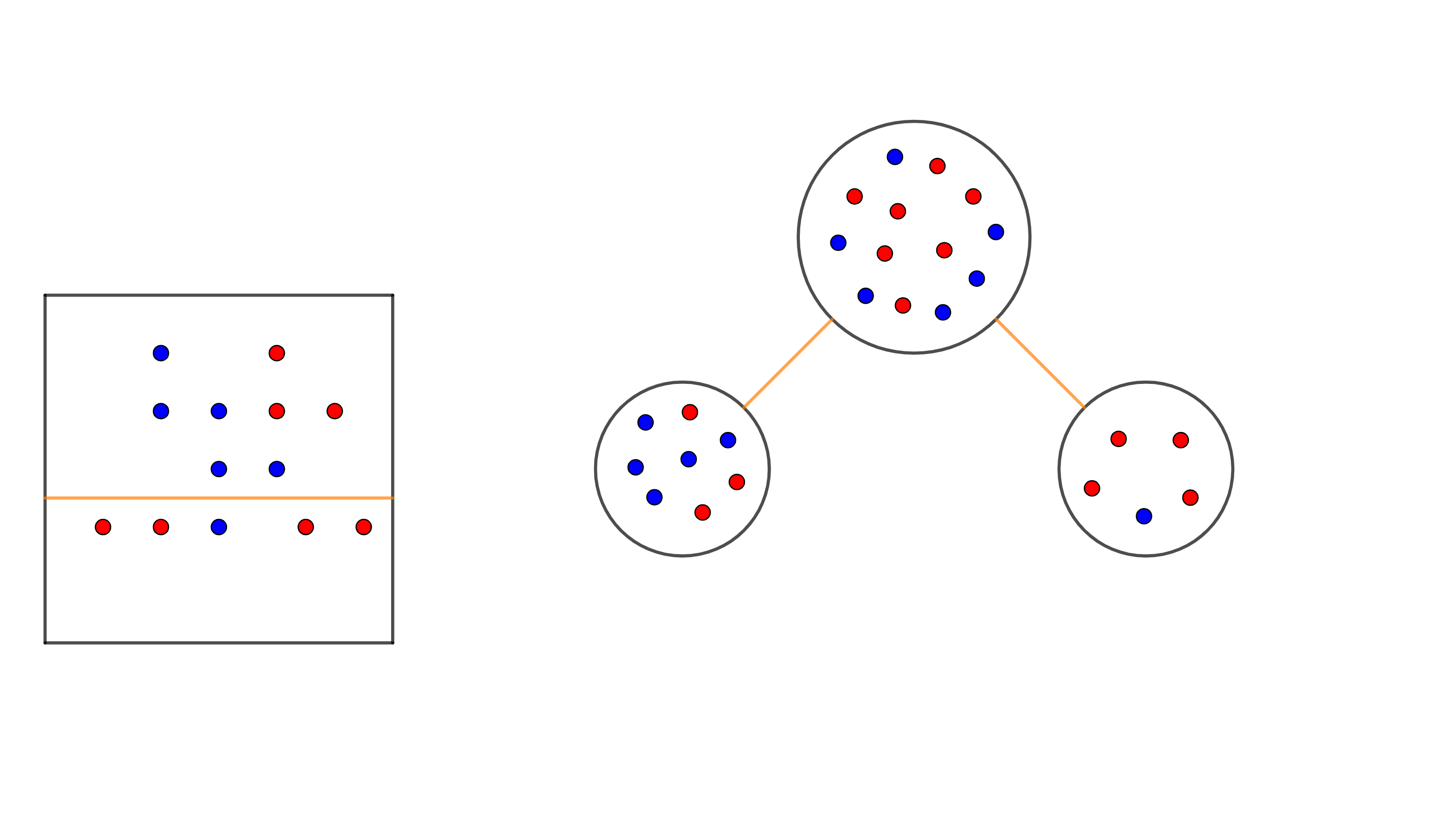}~\includegraphics[width=0.5\textwidth]{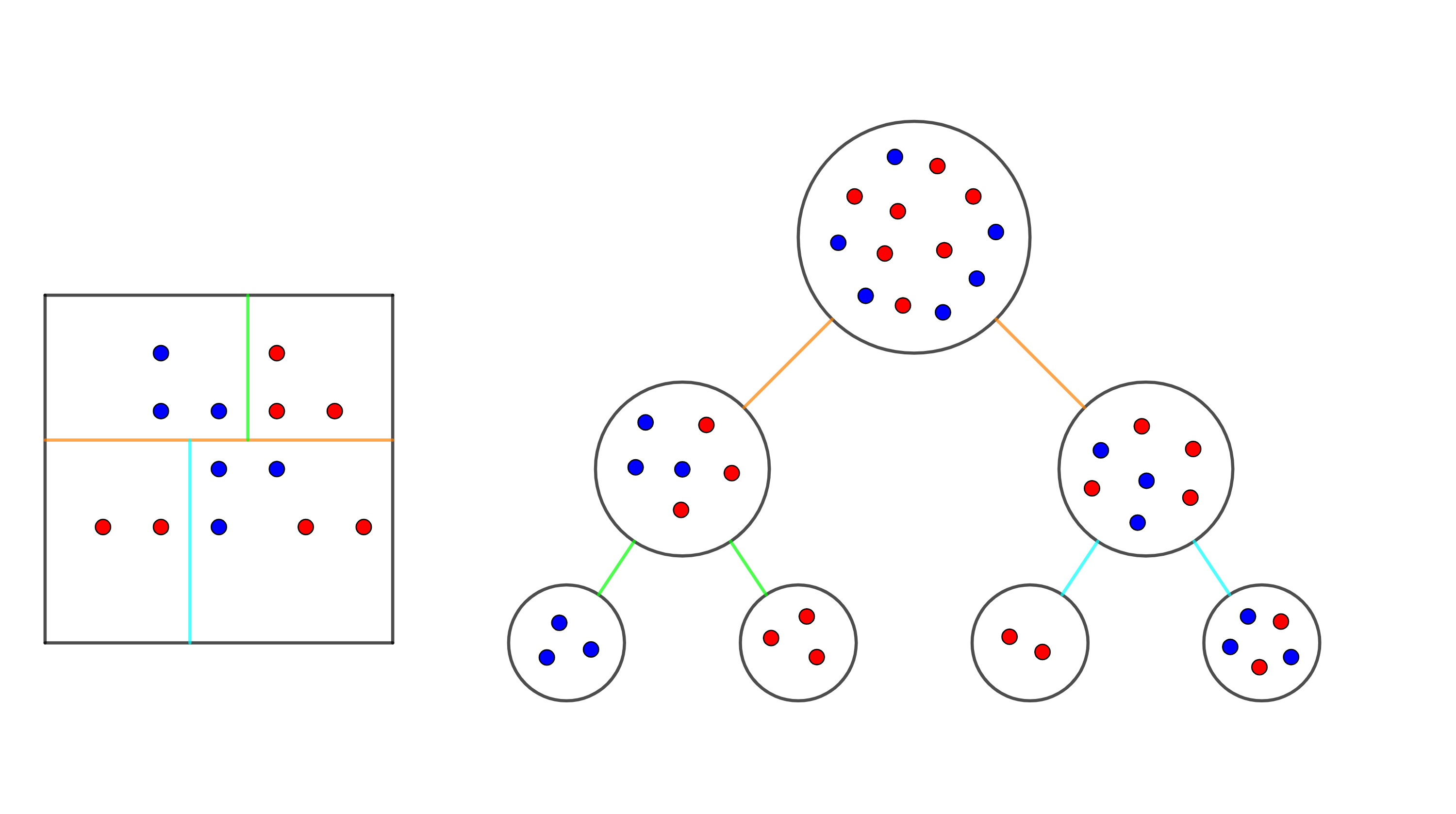}\\
\includegraphics[width=0.5\textwidth]{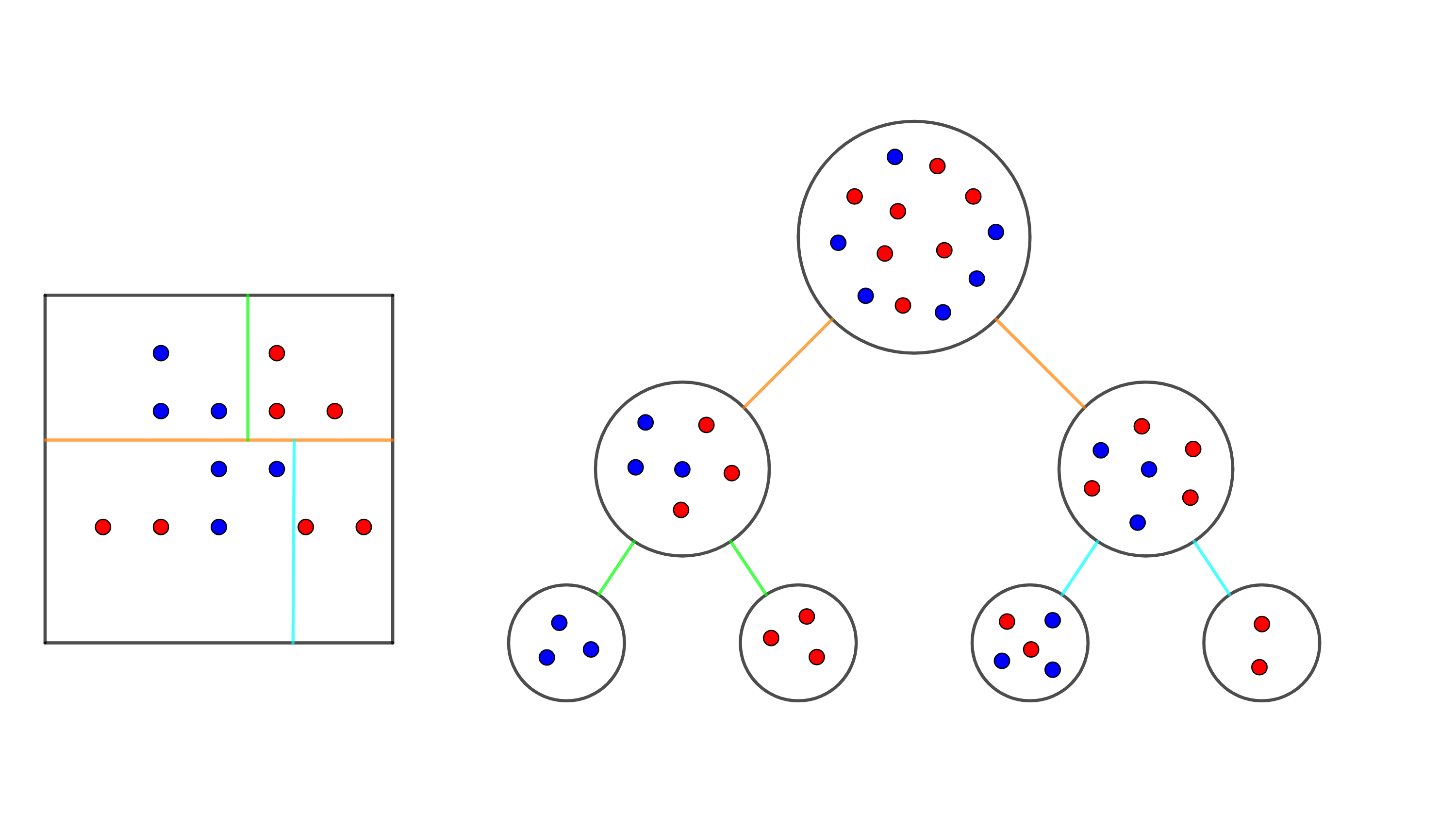}
\caption{Optimal classification forest solution.} \label{fig:2}
\end{figure}

\begin{figure}[h!]
\centering
\includegraphics[width=0.5\textwidth]{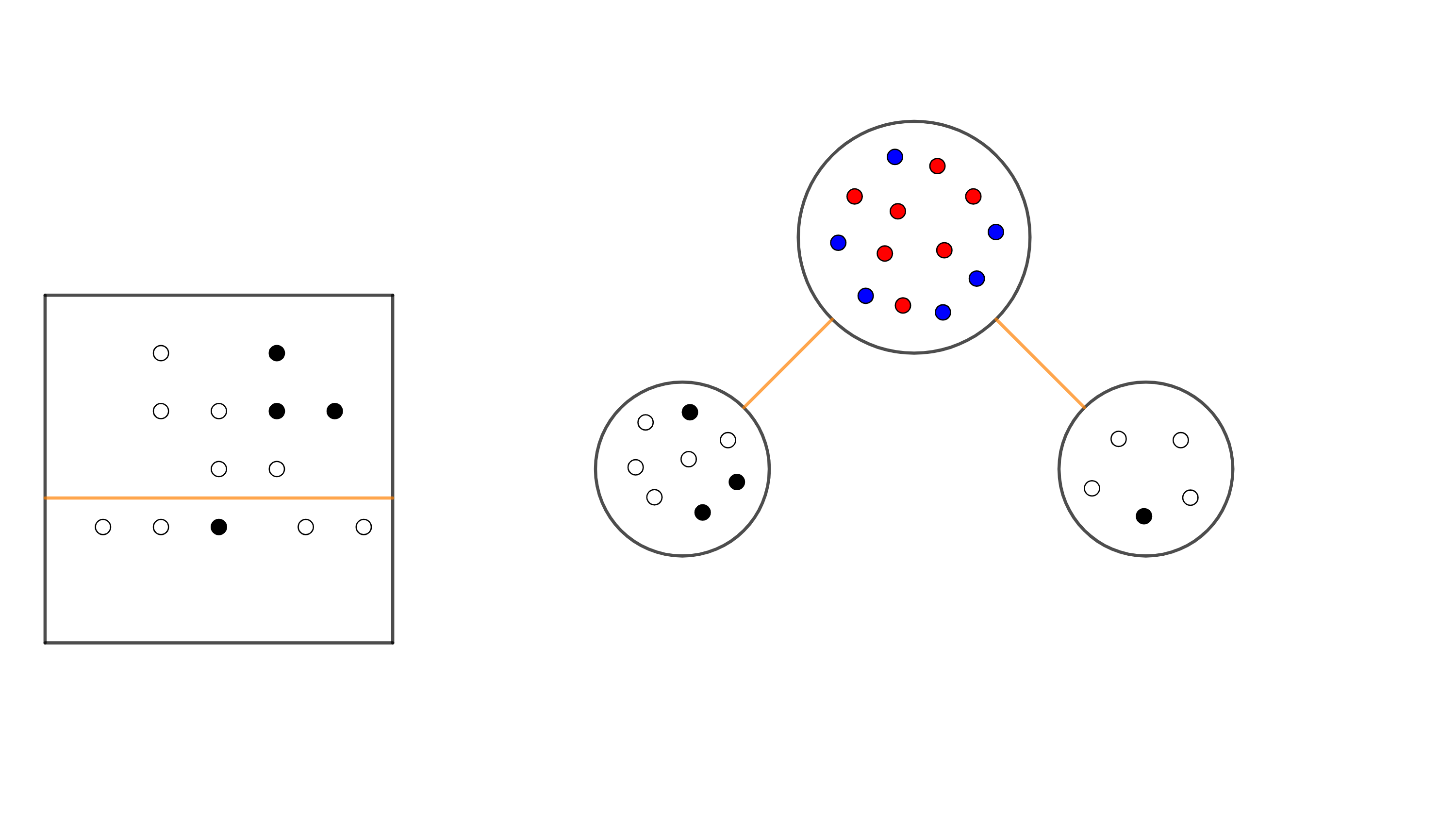}~
\includegraphics[width=0.5\textwidth]{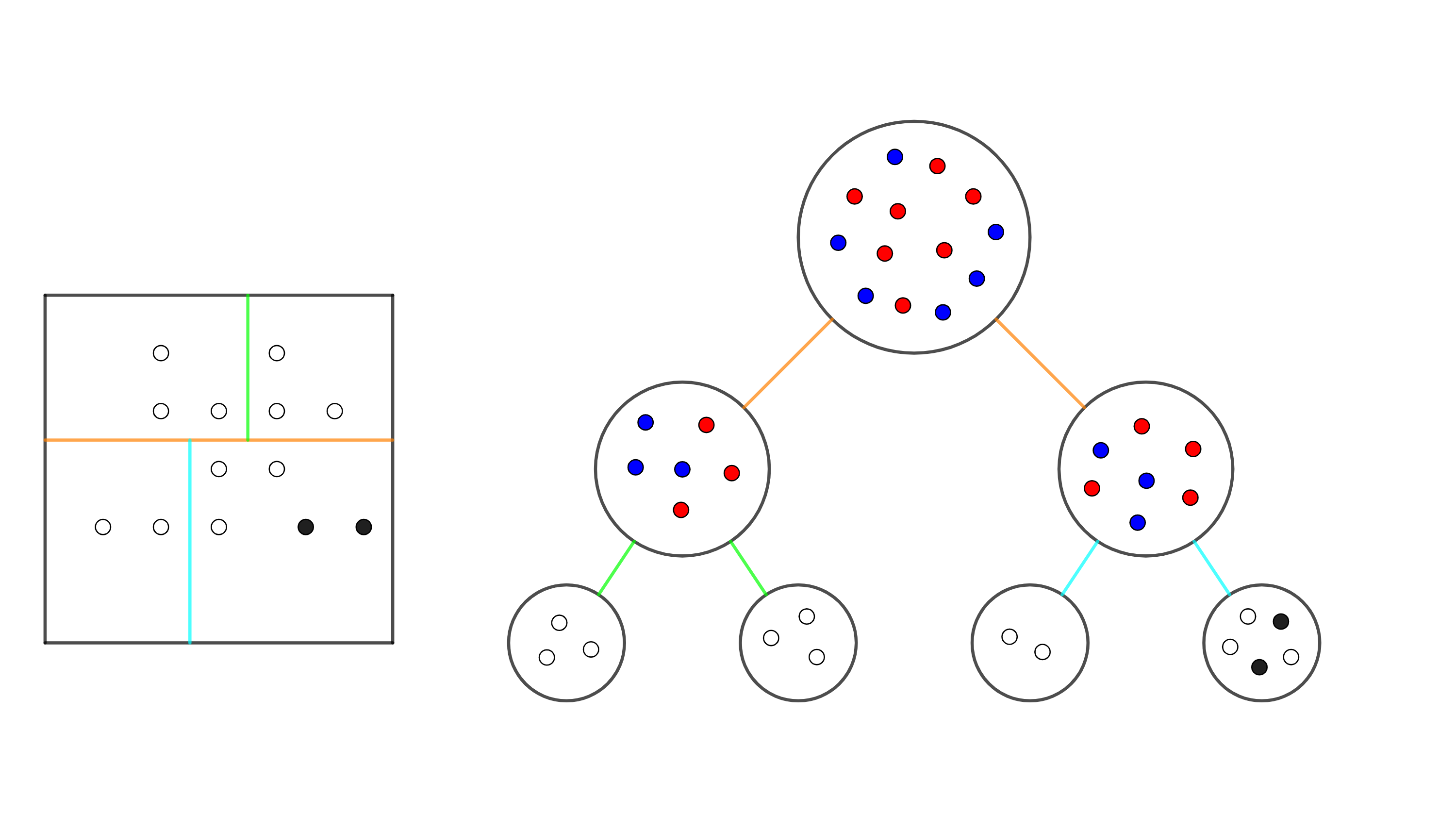}\\
\includegraphics[width=0.5\textwidth]{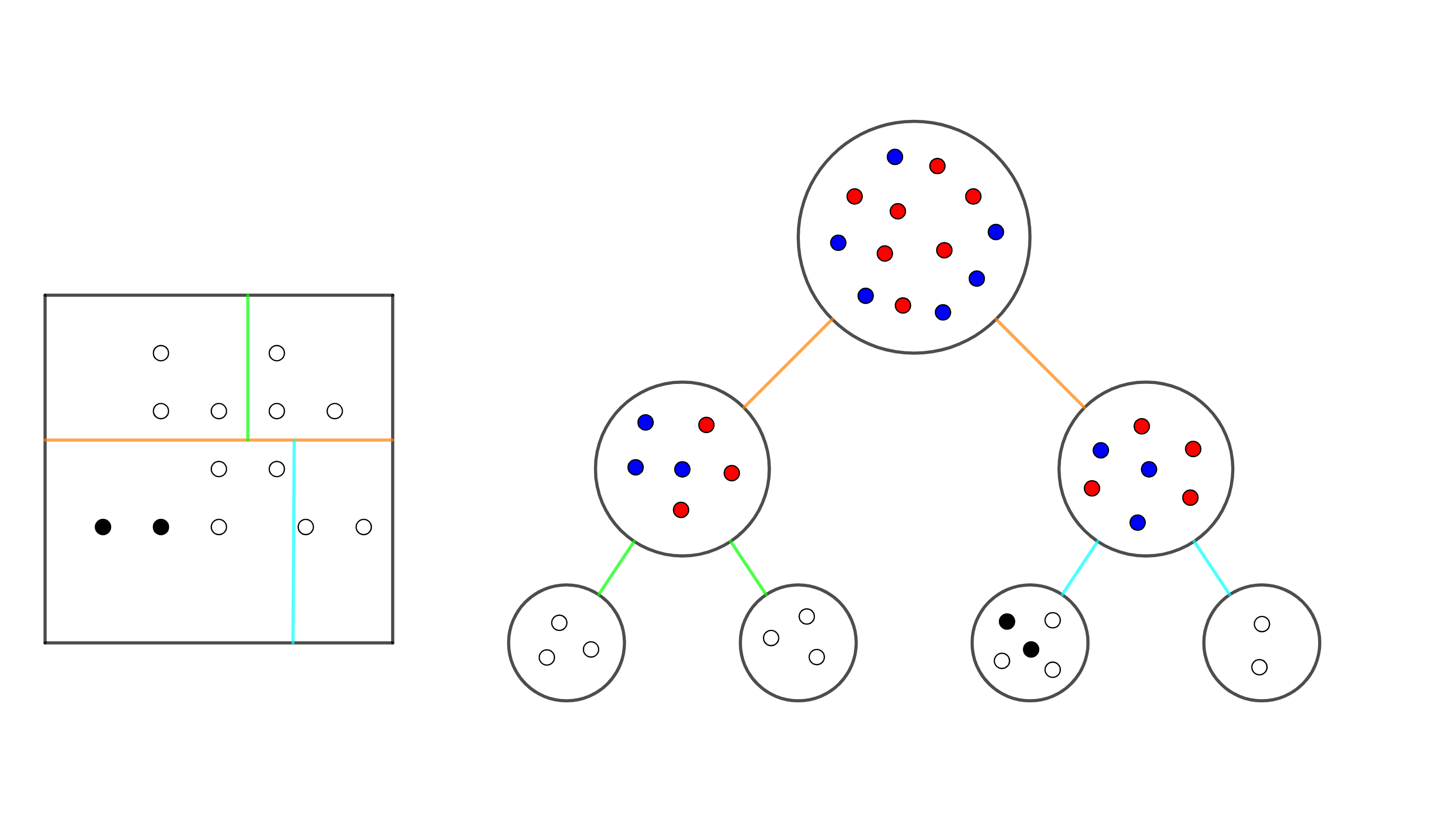}
\caption{Misclassifications occur in individual trees.} \label{fig:3}
\end{figure}

As an example, in Figure \ref{fig:1}  we show a simple OCT with depth two, for a small dataset with $13$ observations. This CT has three branch nodes and four leaf nodes. Besides, it only has one misclassification error, since one blue point is lying in a leaf node where the red class (the most represented class in the leaf) is assigned.
According to this example, to achieve a perfect classification rate, in the context of the OCTs we would have to follow a vertical approach, i.e. increase the OCT depth. In contrast, we propose a horizontal approach through an WOCF to achieve the same end. In this respect, we define an $R$-WOCF (for $R\in \mathbb{Z}_+$) as a set of $R$ CTs that combine their predictions using a (weighted) majority voting strategy in the leaves to construct a global prediction. In this ensemble, all trees have access to the entire training sample and knowledge of what errors are occurring in the other trees in order to coordinate with each other and optimize their strategy. Note that with this definition, the CTs that form an WOCF do not necessarily have to be optimal (in the sense of OCT) for a given depth, and furthermore, neither they have to follow a fixed pattern when it comes to generating splits as in the case of heuristic algorithms, they have to be adequately coordinated to obtain a good final accuracy when voting. To illustrate this paradigm, we show in Figure \ref{fig:2} a solution of the previous example with an WOCF consisting of $3$ CTs of depth $2$.
Looking at each of the trees individually, all of them have non-zero misclassification errors. In Figure \ref{fig:3} we can see the misclassified observations in each tree coloured in black.

Nevertheless, we can see that none of the observations are colored in more than one tree, and therefore, if we were to follow a voting classification strategy, every observation would have at least two votes associated with its correct class and consequently no observation would be finally misclassified.

At this point, one can observe that given a training sample, one can always calibrate the hyperparameters of OCTs and WOCFs to increase their training accuracy. In OCTs, increasing the depth has a direct impact on the training accuracy. In WOCFs, both the depth of the trees and the number of them allow us to increase such accuracy. Splits created in a deep OCT, however, decrease its ``marginal efficiency'', because later splits are conditioned on the previous ones. In other words, deep splits could affect only small subsamples of the training sample, and this could lead to overfitting. On the other hand, the splits in WOCF have less dependency and are more efficient. Thus even if OCT and WOCF have comparable training accuracy, WOCF might achieve better out-of-sample performance. This is indeed the case in the computational experiments we performed, as shown in Section \ref{sec:exps}.

\section{Mathematical Optimization Model}\label{sec:form}

This section is devoted to providing a mathematical optimization formulation of the WOCF for binary instances derived by a forest $\mathcal{F}$ with a fixed odd number of trees, $R$. In order to do this, we start by presenting the objective function of the problem. As it is done in OCT, we minimize the misclassification errors within the training sample as well as the complexity of the ensemble, which is the number of active splits involved in all the trees. Recall that $y_i \in \{0,1\}, \forall i\in N$ are the labels given in the data. We define a binary variable $\alpha_i$ for each observation $i=1,\ldots,n,$ that takes value one in case the ensemble predicts class $1$ for observation $i$, and zero otherwise. Moreover, we use the binary variables $d^r_t \in\{0,1\}$ to model if tree $r$ splits at node $t$. According to these, the objective function of our problem is given by the following expression: 
$$
\min \;\; \frac{1}{n}\sum_{i\in N} \lvert y_i - \alpha_i \rvert + C\sum_{r=1}^{R}\sum_{t\in\tau_b} d_t^r.
$$
That is, the sum of the number of observations that are incorrectly classified (according to the $\alpha$-variables) plus the number of splits in the ensemble. Trees minimizing the above expressions will be constructed.

In what follows we detail the constraints relating to the different variables in our model and that must also ensure that a forest is adequately constructed.

First, we need to ensure that the $\alpha$-variables are properly defined, representing the final result of the weighted voting strategy amongst the $R$ trees. With this goal in mind, we define the binary variables $\theta^r_i, \forall i=1\in N, r=1,\ldots, R,$ which take the value one if class one is predicted in tree $r$ to observation $i$, and takes value zero otherwise, and $w_r \in [0,1]$, for all $r=1, \ldots, R$, that will represent the weight of tree $r$ in the voting.

These variables are adequately defined by means of the following constraints:
\begin{align}
    & \sum_{r=1}^R w_r =1, \label{OCF-0}\\
& \sum_{r=1}^{R} w_r \theta_i^r - \frac{1}{2}  +\epsilon\leq \alpha_i, &  \forall i\in N, \label{OCF-1}\\
& -\sum_{r=1}^{R} w_r \theta_i^r + \frac{1}{2} \leq 1-\alpha_i, &  \forall i\in N,\label{OCF-2}\\
& w_r \geq 0, & \forall r\in\mathcal{F} \label{OCF-2+}.
\end{align}
Constraint \eqref{OCF-0} indicates that the weights for the trees represent proportions. Constraints \eqref{OCF-1} and \eqref{OCF-2} allow the proper definition of the $\alpha$-variable, i.e., if the weighted sum of the votes of all trees to class $1$ for observation $i$  is strictly larger than $50\%$, then, $\alpha_i$ must take value $1$. On the contrary, if the votes are smaller or equal than $50\%$, $\alpha_i$ will take value $0$.

Note that constraints \eqref{OCF-1} and \eqref{OCF-2} are nonlinear, but linearizable by incorporating auxiliary continuous variables, e.g., using McCormick envelopes.

With the above constraints and our objective function, we assure that the voting strategy is well modeled. In this regard, once the basis of the problem has been established, we require each of the trees to be adequately defined. To this end, we  also consider the following set of decision variables to regulate the rationale of the splits in the tree as well as the allocation of the observations to the leaves:

\begin{itemize}
        \item[] $z^r_{it}= \left\{\begin{array}{cl} 1, & \mbox{if observation $i$ belongs to node $t$ in tree $r$,}\\
0, & \mbox{otherwise,}\end{array}\right. $ 
        for all $i \in N$, $t\in\tau_l$, $r\in\mathcal{F}$.
        \item[] $d^r_{t}= \left\{\begin{array}{cl} 1, & \mbox{if a split is applied at node $t$ in tree $r$,}\\
0, & \mbox{otherwise,}\end{array}\right. $ 
for all $t\in\tau_b$, $r\in\mathcal{F}$.
        \item[] $l^r_{t}= \left\{\begin{array}{cl} 1, & \mbox{if $t$ is a non-empty leaf in tree $r$,}\\
0, & \mbox{otherwise,}\end{array}\right. $ 
for all $t\in\tau_l$, $r\in\mathcal{F}$.
    \item[]   $a^r_t \in \{0,1\}^p , b^r_t\in [0,1]$: coefficients and the intercept of the split involved at node $t$ in tree $r$, for all $t\in\tau_b$, $r\in\mathcal{F}$.   
        \item[] $\psi^r_{t}= \left\{\begin{array}{cl} 1, & \mbox{if leaf $t$ of tree $r$ is assigned to class $1$,}\\
0, & \mbox{otherwise,}\end{array}\right. $ 
 \mbox{for all} $ t\in\tau_l$ , $r\in\mathcal{F}$.
\end{itemize}

In contrast to other CTs methodologies track the number of misclassification errors occurred in the leaf nodes of the trees. However, we have to make sure that the observations that belong to the same leaf node in a tree are assigned to the same class. In order to do this, we incorporate to our formulation the following constraints:
\begin{align}
    \psi_t^r - (1-z_{it}^r) \leq \theta_i^r \leq \psi_t^r + (1-z_{it}^r), & \forall i \in N, r\in\mathcal{F}, \ t\in\tau_l.\label{OCF-3}
\end{align}
These constraints assure that in case observations $i$ and $j$ are located in leaf $t$ of tree $r$ ($z_{it}^r=z_{jt}^r=1$), then both $i$ and $j$ would be assigned to class $\psi_t^r$, as $\theta_i^r = \psi_l^r = \theta_j^r  $.

In the following, each of the trees in the ensemble should be adequately defined. We model this structure analogously to other OCT models. For all $r\in \mathcal{F}$:

\begin{align}
& \sum_{q=1}^p a_{tq}^r = d_t^r,  \forall t \in \tau_b, \label{OCF-4}\tag{$6^r$}\\
& 0 \leq b_{t}^r\leq d_t^r,  \forall t \in\tau_b,\label{OCF-5} \tag{$7^r$}\\
& d_t^r \leq d_{p(t)}^r,  \forall\in\tau_b \backslash\{1\},\label{OCF-6}\tag{$8^r$}\\
& \sum_{t\in\tau_l} z_{it}^r=1,  \forall i\in N, \label{OCF-7}\tag{$9^r$}\\
& z_{it}^r \leq l_t^r,  \forall i\in N,\ t\in\tau_l,\ \label{OCF-8}\tag{$10^r$}\\
& \sum_{i=1}^n z_{it}^r \geq n_{\text{min}}l_t^r,  \forall t\in\tau_l,\label{OCF-9}\tag{$11^r$}\\
& (a_m^r)^\intercal x_i + \epsilon \leq b_{m}^r + (1+\epsilon)(1-z_{it}^r), \forall i \in N, \ t\in \tau_l, \label{OCF-10}\tag{$12^r$}\\
& m\in A_L(t),\nonumber\\
& (a_{m}^r)^\intercal x_i \geq b_{m}^r - (1-z_{it}^r),   \forall i\in N, \ t\in \tau_l, \label{OCF-11}\tag{$13^r$}\\
&\ m\in A_R(t). \nonumber
\end{align}

According to this, Constraints \eqref{OCF-4} and \eqref{OCF-5}
enforce that the (univariate) splits are only applied in active nodes $(d^r_t=1)$. Note that oblique splits could be similarly applied by slightly modifying these constraints. Nevertheless,  we maintain orthogonal splits so as to retain the interpretability of the final model. The hierarchical structure of a tree is given in constraint \eqref{OCF-6}. Note that if the parent node of node $t$ does not create a split, then node $t$ should not create a split itself, i.e., $d^r_{p(t)}=0$ implies $d^r_t=0$. In case $t$ creates a split, then its parent node must also create a split,  so $d^r_t=1$ implies  $d^r_{p(t)}=1$. Constraint \eqref{OCF-7} ensures that an observation belongs to exactly one leaf node for each tree. Constraints \eqref{OCF-8} and \eqref{OCF-9} imply that a leaf node is non-empty if and only if there are at least $n_{\text{min}}$ observations in the leaf node. Constraints \eqref{OCF-10}, where  $\epsilon$ is a small positive constant, and \eqref{OCF-11} force observations belonging to leaf $t$ in tree $r$ to satisfy all the splits in order to reach the leaf node $t$ from the root node in tree $r$. Note that in \eqref{OCF-10}, in case $z_{it}^r=0$ (i.e., observation $i$ does not belong to node $t$ in tree $r$), one gets that $(a_m^r)^\intercal x_i  \leq b_{m}^r+1$, which is redundant since $a_m^r \in \{0,1\}^p$, $b_m^{r} \in [0,1]$, and $x_i$ is assumed to be normalized in $[0,1]$. The same reasoning applies to \eqref{OCF-11}.

Gathering all together, the binary WOCF can be formulated as the following MILP problem:

\begin{align*}
\min  \;\; & \;\;\frac{1}{n}\sum_{i=1}^n \lvert y_i - \alpha_i \rvert + C\sum_{r=1}^{R}\sum_{t\in\tau_b} d_t^r  & \tag{${\rm WOCF \ }$}\\
\mbox{s.t.} 
& \ \eqref{OCF-0}-\eqref{OCF-3}, &  \\
&  \eqref{OCF-4}-\eqref{OCF-11},  \forall r \in \mathcal{F}, & \\
& \alpha_i \in\left\{0,1 \right\},  \forall i\in N,\\
& \theta_i^r  \in\left\{0,1 \right\},   \forall i\in N,\ r\in \mathcal{F},\\
& d_t^r \in\left\{0,1 \right\},\ a_t^r \in \left\{0,1 \right\}^p, \ b_{t}^r \in\left[0,1\right],   t\in \tau_b, r\in \mathcal{F},\\
& l_t^r \in\left\{0,1 \right\},   t\in \tau_l, \ r\in \mathcal{F},\\
& z_{it}^r \in\left\{0,1 \right\},  \forall i\in N, \ t \in \tau_l,  r\in \mathcal{F},\\
& \psi_t^r \in \{0,1\}, t\in\tau_l, r\in\mathcal{F}.
\end{align*}

\subsection{Strengthening the model}

The MINLP formulation presented above is valid for our WOCF model. Nevertheless, the problem is computationally expensive to solve, and although it can be handled by most of the off-the-shelf optimization solvers (as Gurobi, CPLEX or XPRESS), the problem can only be solved for small and medium sized instances.

On the other hand, the problem formulation can be extended by a set of valid inequalities that take advantage of intrinsic properties of the problem to help solve it, thus reducing the computational burden. In particular, we have included in our studies some related to the breaking of symmetries in the problem, as it is well-known that symmetric discrete optimization models may be computationally costly since deeper branch-and-bound trees are to be
constructed to guarantee the optimality of the solutions.

Note that for a given classification problem, each feasible forest is equivalent (in objective value) to all permutations in the tree indices $r$. The construction of the trees is in a sense free in the ensemble since the indexing in $r$ is only useful to enumerate the trees. However, without affecting the correct functionality of the model, these trees could be sorted by classification errors, so that the first index in $r$ coincides with the most accurate tree and the last index in $r$ with the least. In this way, we would avoid permutations of trees being valid solutions of the problem, and thus breaking the symmetries of the problem. This ordering can be done by the following constraints:
\begin{align*}
& \displaystyle\sum_{i=1}^n \left| y_i- \theta_{i}^1 \right| \leq \displaystyle\sum_{i=1}^n \left| y_i-\theta_{i}^2 \right| \leq \ldots \leq \displaystyle\sum_{i=1}^n \left| y_i-\theta_{i}^R \right|.
\end{align*}
The solvers mentioned above do not support this type of constraints directly, and therefore the absolute value constraints must be replaced, with the help of some auxiliary binary variables, $\lambda_{i}^r \in \left\lbrace 0,1 \right\rbrace$, by the following sets of constraints:
\begin{align}
&  y_i-\theta_i^r \leq \lambda_i^r, \label{sym_1}\\
&  -y_i + \theta_i^r \leq \lambda_i^r, \label{sym_2}\\
& \lambda_i^r \leq y_i+\theta_i^r  ,\label{sym_3}\\
& \lambda_i^r \leq 2 - y_i -\theta_i^r .\label{sym_4}
\end{align}

Thus, given that $\theta_i^r, y_i \in \left\lbrace 0,1 \right\rbrace $, the $\lambda_i^r$ variables  represent exactly the absolute value above. Note that if $\theta_i^r\neq y_i$, either Constraints \eqref{sym_1} or \eqref{sym_2} will be activated forcing $\lambda_i^r$ to take value $1$. Besides, if $\theta_i^r = y_i = 0$, Constraints \eqref{sym_3} will force $\lambda_i^r$ to take value $0$, and the same will apply to Constraints \eqref{sym_4} if $\theta_i^r = y_i = 1$. Thus, we finally manage to eliminate the symmetries of the permutations of the trees by adding the following constraints on the problem: 
\begin{align*}
& \displaystyle\sum_{i=1}^n\lambda_i^1 \leq \displaystyle\sum_{i=1}^n \lambda_i^2  \leq \ldots \leq \displaystyle\sum_{i=1}^n  \lambda_i^r.
\end{align*}

The above symmetry breaking constraints are not the only constraints that one can apply to reduce the search in the branch-and-bound tree when solving the problem. Instead of sorting the labels given to the trees by the classification errors, one could also sort by the largest error obtained in any leaf of the tree, by the number of splits applied in each tree, or by the largest/smallest intercept ($\max_{t} b_{t}^r$) applied. All these alternatives can be formulated as linear constraints using auxiliary variables as the one detailed above. Nevertheless, after testing some preliminary instances to select the best strategy, we concluded that their effect was similar.

On the other hand, note that in our model, the constraints that we introduce to regulate the logistics of a CT, \eqref{OCF-4}-\eqref{OCF-11}, are inherited from \citep{OCT}, where the first general formulation of an OCT was presented. However, alternative formulations have been recently provided to construct OCTs.  One may replace the above mentioned block of constraints (for each tree in the ensemble) by any of the available alternatives. Specifically, \cite{verwer2019learning} propose \texttt{binOCT*}, a binary linear integer optimization formulation which avoids using the whole set of feature values but uses a binary encoding of the different thresholds induced by the training dataset. \cite{flowOCT} provide a maxflow based formulation (\texttt{flowOCT}), together with a Benders decomposition algorithm, for datasets with binary features. The formulation  proposed in \citep{gunlunk} also exploits the combinatorial nature of datasets with categorical features. With the same paradigm that OCT, several extensions have been also proposed (see e.g., \cite{strong,MOCTSVM,OCTSVM,stable_class,fair,zhu}, among many others).

Nevertheless, the formulation in \cite{OCT} has been proven to produce a reasonable trade-off between high performance, versatility (in the types of datasets that can be classified), and training CPU time. Thus, we decided to present our model under this framework, although it can be adapted to a different formulation, in case the dataset meets the assumptions required for its use.

\section{Computational Experiments}\label{sec:exps}

In this section, we present the results of our computational experiments, with the aim of studying the predictive capacity of WOCF with respect to other similar methods. Since the idea is to obtain interpretable classifiers, as we show in the next section, we study simple versions of the possible resulting models. Nine methods are compared: CART, OCT, RF with 3 trees (3-RF) and 500 trees (500-RF), XGBoost with 3 trees (3-XGB), 5 trees (5-XGB) and 500 trees (500-XGB), and WOCF with 3 trees (3-WOCF) and 5 trees (5-WOCF). All these methods have been applied to  twenty datasets from the UCI Machine Learning Repository~\citep{UCI}. The datasets' names, abbreviations, and dimensions ($n$: number of observations, $p$: number of features) are reported in Table \ref{table:0}. We note that these datasets are neither novel nor impressive in size or complexity, sometimes even considered too well-tested. But that is precisely the point of this section -- we want to demonstrate the performance of WOCF on these well-tested datasets and show from ground up that WOCFs can provide Pareto improvement over previous methods. In the next section, we will study WOCFs on more complex datasets and provide more insights. 

\begin{table}[h!]
\centering
\begin{tabular}{@{}lcc@{}}\hline
Dataset & Abbreviation & (Samples, Features) \\\hline
\texttt{Connectionist Bench (Sonar, Mines vs. Rocks)} & sonar &(208, 60) \\
\texttt{Parkinson Dataset with replicated acoustic features} & parkinson &(240, 40) \\
\texttt{Statlog (Heart)} & heart & (270, 13) \\
\texttt{Vertebral Column} & vertebral & (310, 6) \\
\texttt{Ionosphere} & ionosphere &(351, 33) \\
\texttt{Thoracic Surgery Data} & thoracic &(470, 16) \\
\texttt{Early Stage Diabetes Risk Prediction} & diabetes &(520, 16) \\
\texttt{ILPD (Indian Liver Patient Dataset)} & ildp & (552, 10) \\
\texttt{Breast Cancer Wisconsin (Original)} & bcw & (683, 9) \\
\texttt{SPECT Heart} & s-heart& (683, 9) \\
\texttt{Statlog (australian Credit Approval)} & australian & (690, 14) \\
\texttt{Mammographic Mass } & mmass &(830, 6) \\
\texttt{Glioma} & glioma & (862, 21) \\
\texttt{Raisin} & raisin & (900, 7) \\
\texttt{Tic-Tac-Toe} & ttt & (958, 8) \\
\texttt{Statlog (German Credit Data)} & german &(1000, 24) \\
\texttt{QSAR biodegradation} & biodeg &(1055, 41) \\
\texttt{Auction verification} &  auction & (2048, 7) \\
\texttt{Chess (King-Rook vs. King-Pawn)} & chess & (3195, 37) \\
\texttt{Rice (Cammeo and Osmancik)} & rice & (3810, 7) \\\hline
\end{tabular}
\caption{Dataset details.}
\label{table:0}
\end{table}

To establish a fair comparison over the datasets where all observations are involved at some point as a training, validation or test set, we perform a 4-fold cross-validation scheme, i.e., datasets are split into four random train-validation-test partitions. Two of the folds are used for training the model, one of them is used as a validation set for calibrating the parameters of the models, and the remaining one is used for testing. For each method, we choose the combination of parameters with a better performance (in terms of accuracy) in the validation sample and use that model to compute the accuracy, in percentage, for the test set.  Finally, in order to avoid taking advantage of beneficial initial  train-validation-test partitions, we repeat the cross-validation process 5 times and report the average accuracy and standard deviation.\\

With regard to the parameters involved in CART and OCT, we set the minimum amount of observations allowed to be in a leaf node, $n_{\text{min}}$, to be $2.5\%$ of the training sample size, and the depth of the trees, $D$, is set to three. Moreover, as it is done by \citet{OCT}, instead of calibrating a complexity parameter, we remove the complexity term from the objective function and add the following constraint that upper bounds the number of active splits:
\begin{align*}\sum_{t\in\tau_b} d_t \leq C.
\end{align*}
where $C \in \{1,\ldots,9\}$ is the parameter to be tuned during validation.

For the 3-WOCF and 5-WOCF, we use a smaller depth, $D = 2$, and set $n_{\text{min}}$ to be $2.5\%$ of the training sample. We let the \textit{total} number of splits to be in $C \in \{3,\ldots,9\}$ for the 3-WOCF, and $C \in \{5,7,\ldots,15\}$ for the 5-WOCF, and we regulate the size of the forest by establishing: 
\begin{align*} \sum_{r =1}^{R} \sum_{t\in\tau_b} d_t^r \leq C.
\end{align*}
The 3-RF and 500-RF models were trained using the \texttt{scikit-learn} library, where the hyperparameters regulating the forest size were applied as in the WOCF case, and the rest were left as default. The same applies for the 3-XGB, 5-XGB and 500-XGB models, which were trained with the \texttt{XGBoost} library. In the latter, no $L_1$ or $L_2$ regularization was applied on the weights to equalize conditions with the WOCF models.\\

Additionally, both OCT and WOCF formulations have a family of constraints that involves a small enough parameter, which has been fixed to the value $\epsilon=10^{-5}$ in all the experiments. 
These models were coded in \texttt{Python} and solved using the optimization solver \texttt{Gurobi v11.0.1}.  The experiments were run on a PC Intel Xeon E-2146G processor at 3.50GHz and 64GB of RAM. In all cases, a time limit of one hour was set for training and in most of the instances this limit was reached without guaranteeing optimality.\\

We present two tables with the results of all experiments, where we report the average accuracy (and standard deviations) of all methods for each of the twenty datasets. In Table \ref{table:1}, we show the comparison between 3-WOCF, CART, OCT, 3-RF and 500-RF. In table \ref{table:2}, we present a more extensive comparison of the ensemble models with weights in the trees, covering the models 3-WOCF, 3-XGB, 5-WOCF, 5-XGB and 500-XGB. In the latter, the results are boldfaced in green when a model, WOCF or XGB, performs better by more than $1\%$ than its competitor for the same number of trees.
\\
The most remarkable observation, in view of the results, is that 3-WOCF tends to outperform the other classification tree methods in most of the datasets. Moreover, it achieves results comparable to those of 500-RF. In this sense, in these small datasets, we show a Pareto improvement in terms of accuracy versus interpretability. We obtain better results in accuracy compared to the interpretable models (3-RF and CTs) and better results in terms of interpretability compared with complex random forests, as 500-R.\\
On the other hand, as for Table \ref{table:2}, we observe that in general the 3-WOCF model performs better than the 3-XGB model, obtaining in seven of the twenty datasets an accuracy higher than $1\%$ than that of XGB, and only being below in two of them. However, when dealing with 5 trees, although the results are fairly close, in five datasets the performance of XGB is better than that of WOCF, and only in two of them its performance is worse. Note that some of these data sets are very small, and these are cases where large MIP models may show more overfitting than other methods. In fact, in three of the datasets where 5-XGB outperforms 5-WOCF, the 500-XGB model also underperforms the other models, suggesting that this data may have been more susceptible to overfitting. Furthermore, in larger datasets, as optimality is not guaranteed, and considering the fact that the 5-WOCF model is more difficult to solve than the 3-WOCF, solutions may have room for improvement with longer time limits. Finally, it should also be noted that the 500-XGB model results in performances very difficult to match, and although the dominance of this model is not total, the WOCF model with low number of trees has not been able to match 500-XGB model in terms of accuracy, as it does when compared to the RF model. However, explainability and interpretability of WOCFs are, by construction, superior to other tree ensemble methods as RF or XGB.

\begin{table}[h!]
\centering
\begin{adjustbox}{width=\textwidth}
\begin{tabular}{@{}l*{6}{c}@{}}\hline
Dataset & CART & OCT & 3-RF & 500-RF & 3-WOCF \\\hline
\texttt{sonar} & 65.02 $\pm$ 5.90 & 71.04 $\pm$ 7.06 & 69.10 $\pm$ 9.20 & 72.72 $\pm$ 6.39 & 73.48 $\pm$ 8.23  \\
\texttt{parkinson} & 72.13 $\pm$ 6.66 & 75.19 $\pm$ 5.97 & 75.94 $\pm$ 6.60 & 77.72 $\pm$ 5.09 & 76.04 $\pm$ 6.35  \\ 
\texttt{heart} & 74.36 $\pm$ 4.87 & 76.91 $\pm$ 4.49 & 74.45 $\pm$ 5.22 & 80.16 $\pm$ 4.94 & 77.53 $\pm$ 5.90  \\
\texttt{vertebral} & 79.75 $\pm$ 4.43 & 80.13 $\pm$ 4.35 & 81.15 $\pm$ 4.38 & 81.94 $\pm$ 4.42 & 80.65 $\pm$ 4.99  \\
\texttt{ionosphere} & 88.11 $\pm$ 3.28 & 88.55 $\pm$ 3.66 & 88.52 $\pm$ 4.04 & 90.68 $\pm$ 3.83 & 90.40 $\pm$ 3.36 \\
\texttt{thoraric} & 82.85 $\pm$ 3.15 & 83.75 $\pm$ 2.87 & 84.76 $\pm$ 2.99 & 85.15 $\pm$ 2.71 & 83.95 $\pm$ 2.60 \\
\texttt{diabetes} & 89.07 $\pm$ 3.17 & 91.34 $\pm$ 2.32 & 87.10 $\pm$ 3.39 & 89.15 $\pm$ 2.69 & 90.92 $\pm$ 2.69  \\
\texttt{ildp} & 67.19 $\pm$ 3.93 & 68.15 $\pm$ 2.80 & 67.31 $\pm$ 3.53 & 68.25 $\pm$ 3.59 & 66.76 $\pm$ 4.14 \\
\texttt{bcw} & 92.14 $\pm$ 1.81 & 94.24 $\pm$ 1.89 & 94.28 $\pm$ 2.21 & 96.21 $\pm$ 1.50 & 95.18 $\pm$ 1.77  \\
\texttt{s-heart} & 75.79 $\pm$ 5.38 & 76.69 $\pm$ 4.57 & 76.40 $\pm$ 6.19 & 78.27 $\pm$ 5.57 & 78.66 $\pm$ 5.28  \\
\texttt{australian} & 84.16 $\pm$ 1.71 & 84.36 $\pm$ 2.28 & 80.06 $\pm$ 4.83 & 85.12 $\pm$ 2.10 & 85.39 $\pm$ 2.46  \\
\texttt{mmass} & 82.98 $\pm$ 2.30 & 82.70 $\pm$ 2.01 & 80.40 $\pm$ 2.30 & 84.13 $\pm$ 2.58 & 83.24 $\pm$ 2.07  \\
\texttt{glioma} & 86.85 $\pm$ 1.94 & 86.24 $\pm$ 2.08 & 84.38 $\pm$ 1.87 & 86.21 $\pm$ 2.11 & 87.03 $\pm$ 1.87  \\
\texttt{raisin} & 85.01 $\pm$ 2.29 & 85.14 $\pm$ 1.77 & 85.38 $\pm$ 2.26 & 85.94 $\pm$ 2.10 & 85.28 $\pm$ 1.91  \\
\texttt{ttt} & 69.93 $\pm$ 2.92 & 72.77 $\pm$ 1.97 & 69.08 $\pm$ 3.17 & 67.32 $\pm$ 3.62 & 75.99 $\pm$ 2.96 \\
\texttt{german} & 70.70 $\pm$ 3.14 & 71.42 $\pm$ 3.32 & 69.65 $\pm$ 3.79 & 72.07 $\pm$ 2.27 & 71.06 $\pm$ 2.61  \\
\texttt{biodeg} & 78.76 $\pm$ 2.59 & 78.21 $\pm$ 2.28 & 77.74 $\pm$ 2.76 & 79.29 $\pm$ 2.40 & 78.13 $\pm$ 2.90 \\
\texttt{auction} & 91.11 $\pm$ 1.46 & 90.85 $\pm$ 1.88 & 87.54 $\pm$ 1.07 & 87.16 $\pm$ 1.01 & 91.05 $\pm$ 1.53  \\
\texttt{chess} & 90.46 $\pm$ 1.04 & 90.68 $\pm$ 1.86 & 90.21 $\pm$ 1.97 & 93.55 $\pm$ 1.14 & 92.11 $\pm$ 2.99  \\
\texttt{rice} & 92.52 $\pm$ 0.70 & 92.53 $\pm$ 0.67 & 92.78 $\pm$ 0.70 & 92.43 $\pm$ 0.52 & 92.61 $\pm$ 0.74  \\\hline
\end{tabular}
\end{adjustbox}
\caption{Average accuracy $\pm$ std obtained in our computational experiments for CART, OCT, RF and 3-WOCF.}\label{table:1}
\end{table}

\begin{table}[h!]
\centering
\begin{adjustbox}{width=\textwidth}
\begin{tabular}{@{}l*{6}{c}@{}}\hline
Dataset & 3-WOCF & 3-XGB & 5-WOCF & 5-XGB & 500-XGB \\\hline
\texttt{sonar} &  \color{armygreen}\bf 73.48 $\pm$ 8.23 & 69.20 $\pm$ 2.35 &
 69.23 $\pm$ 8.45 & \color{armygreen}\bf 70.47 $\pm$ 7.59 & 80.11 $\pm$ 5.94 \\
\texttt{parkinson} & 76.04 $\pm$ 6.35 & 75.25 $\pm$ 6.52 &  76.53 $\pm$ 4.95 & 76.70 $\pm$ 6.67 & 76.70 $\pm$ 5.02 \\ 
\texttt{heart} &  77.53 $\pm$ 5.90 & 78.14 $\pm$ 5.58 & 78.01 $\pm$ 4.83 & \color{armygreen}\bf 79.55 $\pm$ 5.73 & 79.56 $\pm$ 4.75\\
\texttt{vertebral} & \color{armygreen}\bf 80.65 $\pm$ 4.99 & 78.84 $\pm$ 4.82 &
80.71 $\pm$ 3.98 & 80.46 $\pm$ 4.83 & 80.71 $\pm$ 4.78 \\
\texttt{ionosphere} & \color{armygreen}\bf 90.40 $\pm$ 3.36 & 88.80 $\pm$ 3.40 &
89.42 $\pm$ 4.23 & 90.17 $\pm$ 3.10 & 91.77 $\pm$ 2.07\\
\texttt{thoraric} & 83.95 $\pm$ 2.60 & \color{armygreen}\bf 85.06 $\pm$ 2.68 &
83.18 $\pm$ 2.81 & \color{armygreen}\bf 84.59 $\pm$ 4.23  & 79.70 $\pm$ 3.15 \\
\texttt{diabetes} & \color{armygreen}\bf 90.92 $\pm$ 2.69 & 88.57 $\pm$ 2.10 &
\color{armygreen}\bf 92.51 $\pm$ 2.14 & 88.53 $\pm$ 1.80  & 94.48 $\pm$ 2.69 \\
\texttt{ildp} & 66.76 $\pm$ 4.14 & \color{armygreen}\bf 69.38 $\pm$ 3.40 &
67.85 $\pm$ 3.63 &  \color{armygreen}\bf 68.94 $\pm$ 2.92  & 67.53 $\pm$ 3.87\\
\texttt{bcw} &  95.18 $\pm$ 1.77 & 94.51 $\pm$ 1.81 &
95.33 $\pm$ 1.68 & 95.39 $\pm$ 1.97 &
96.01 $\pm$ 1.33\\
\texttt{s-heart} & 78.66 $\pm$ 5.28 & 79.11 $\pm$ 4.25 &
81.27 $\pm$ 3.68 & 80.86 $\pm$ 5.76 &
78.88 $\pm$ 4.77 \\
\texttt{australian} & 85.39 $\pm$ 2.46 & 85.44 $\pm$ 2.35 &
85.16 $\pm$ 2.41 & 85.21 $\pm$ 2.41 &
85.04 $\pm$ 1.77 \\
\texttt{mmass} & 83.24 $\pm$ 2.07 & 83.89 $\pm$ 2.69 &
82.57 $\pm$ 1.98 & \color{armygreen}\bf 83.68 $\pm$ 2.61 & 80.18 $\pm$ 1.80\\
\texttt{glioma} & 87.03 $\pm$ 1.87 & 86.91 $\pm$ 1.99 &
86.54 $\pm$ 2.41 & 86.98 $\pm$ 2.01 &
85.24 $\pm$ 2.21\\
\texttt{raisin} & 85.28 $\pm$ 1.91 & 85.30 $\pm$ 2.40 &
84.56 $\pm$ 2.23 & 85.76 $\pm$ 2.09 &
83.63 $\pm$ 1.80 \\
\texttt{ttt} & \color{armygreen}\bf 75.99 $\pm$ 2.96 & 69.43 $\pm$ 3.02 &
\color{armygreen}\bf 77.76 $\pm$ 3.72 & 71.93 $\pm$ 2.98 & 91.19 $\pm$ 1.50\\
\texttt{german} & 71.06 $\pm$ 2.61 & 70.72 $\pm$ 2.50 &
71.82 $\pm$ 2.53 & 71.72 $\pm$ 2.53 &
73.70 $\pm$ 2.04\\
\texttt{biodeg} & 78.13 $\pm$ 2.90 & 78.72 $\pm$ 2.20 &
80.10 $\pm$ 2.53 & 81.02 $\pm$ 2.13 &
84.97 $\pm$ 1.59 \\
\texttt{auction} & \color{armygreen}\bf 91.05 $\pm$ 1.53 & 88.54 $\pm$ 1.75 &
89.94 $\pm$ 1.93 & 89.96 $\pm$ 1.27 & 
98.84 $\pm$ 0.82\\
\texttt{chess} & \color{armygreen}\bf  92.11 $\pm$ 2.99 & 90.79 $\pm$ 1.55 &
92.95 $\pm$ 1.90 & 93.66 $\pm$ 1.28 &
98.79 $\pm$ 0.44 \\
\texttt{rice} & 92.61 $\pm$ 0.74 & 92.72 $\pm$ 0.76 & 92.72 $\pm$ 0.83 & 92.76 $\pm$ 0.75 & 91.23 $\pm$ 0.95\\
\hline
\end{tabular}
\end{adjustbox}
\caption{Average accuracy $\pm$ std obtained in our computational experiments for methods WOCF and XGB.}\label{table:2}
\end{table}

\section{Case Studies}\label{sec:case}

In the previous section we studied the predictive ability of our model with respect to its tree-shaped benchmarks. We showed how the model is able to broadly improve the accuracy of interpretable CT-based models, but we did not discuss in detail the obtained solutions. As we have pointed out throughout this paper, in the context of machine learning, solutions that are interpretable by humans are preferable to black-box approaches where it is not possible to identify why a certain class was predicted to a given observation. \\

In practice we often find that for certain datasets numerous classifiers report very similar accuracy results (this is indeed the case in some of our computational experiments), and in these cases there is no dispute that it is better to use the simpler models. On the other hand, the work presented by \cite{IML} points out that when several classifiers of different nature report the same accuracy, there should be a simple and interpretable model that preserves accuracy. However, finding such a simple model, in case of existence, can be very challenging. \\

In this section we present the analysis of three real case studies. We show in detail the solutions obtained with the WOCFs and CTs and make a comparative analysis of them. One of the advantages of WOCF, apart from the increase in accuracy, is the flexibility when interpreting the solutions, as the voting system is less rigid than following the paths through the branches of a single tree. This fact is due to the independence of the variable clusters on different trees of the ensemble. As a result, we show that in some cases WOCF can provide some of the simplest and most practical models to explain the rationale of a classification rule for a given data set, understanding what the important factors are and how they are grouped together.\\

The major disadvantage of the WOCF model is its computational cost, due to the difficulty of solving a MILP where the number of binary variables increases as the size of the dataset increases. In the previous section, we reported a symmetry-breaking constraint to overcome this problem. In order to make a fair comparison, the technique used was a simple random sampling, however, this technique can be refined to efficiently target the resolution of large datasets by solving optimization problems involving a few sets of points.\\

For the sake of simplicity in the interpretations in these case studies, we fixed the values for the weights $w_r$ all equal, specifically $w_r=\frac{1}{3}$, for all $r=1, 2, 3$.

\subsection{Case Study 1: Granted Mortgages}

In this first study, we analyze a problem presented by Dream Home Financing, which is  a company that provides loans to clients to enable them to finance their houses. The objective of the problem is to find a classification rule to determine whether a loan is granted to a client or not. The dataset, which is freely available, can be found on the Analytics Vidhya competition (\href{https://datahack.analyticsvidhya.com/contest/practice-problem-loan-prediction-iii/#About}{\texttt{grant\_mortgages}}). The dataset, after filtering out the rows containing any N/A among the variables, collects the information concerning 480 clients with a distribution of the target variable of $69\%$ (value 1 = loan granted) and $31\%$ (value 0 = loan not granted). Besides, it contains the following $11$ predictive variables:

\begin{table}[h!]
\centering
\begin{tabular}{ll}
    Feature & Description\\\hline
   \texttt{Gender} & Male / Female.\\
   \texttt{Married} & Yes / No.\\
   \texttt{Dependents} & Number of dependents.\\
   \texttt{Education} & Graduated / Not Graduated.\\
   \texttt{Self Employed} & Yes / No.\\
   \texttt{Applicant Income} & Continuous variable.\\
   \texttt{Co-applicant Income (Co-Income)} & Continuous variable.\\
   \texttt{Loan Amount (L.A.)} & Continuous variable.\\
   \texttt{Loan Amoun Term (L.A.T.)} & Continuous variable.\\
   \texttt{Credit History} & Yes / No (credit history meets the guidelines).\\
   \texttt{Property type} & Rural / Semi-Urban / Urban.\\\hline
\end{tabular}
\caption{Predictive variables for the Grant Mortgage Dataset.\label{vars:mortgage}}
\end{table}
Note that both parties have an interest in the loan being granted, the client for their personal benefit and the company for the economic benefit it generates. Now, if we look at the set of variables, we see that many of them have values that are difficult for the client to change instantaneously, such as the annual income or the level of education. Nevertheless, the variable L.A.T. is easily modifiable up to a certain upper threshold by the company, and even L.A. can have small modifications without compromising the viability of the operation on the part of the client. \\

On the other hand, if we consider clients who are not granted the loan, we know that by modifying certain values in their input variables the loan could have been granted. However, these changes may be due to variables over which we have room for maneuver or over which we do not. Once a model has been trained, if it is a black-box model, the fact that a variable can be modified according to the business is irrelevant, but on the other hand, if the model is interpretable and such a model involves variables that can be modified, the company could rescue some of the clients who were not originally granted the loan and thus increase its profits (consider, for example, a client who is originally refused a loan of $X$ amount in time $T$, but who could be granted the loan in time $T+t$). As a result of an interpretable model, one should be able  to derive counterfactual explanations to clients which are unsatisfied with their classification.

\begin{figure}[h!]
\centering
\fbox{\includegraphics[width=0.8\textwidth]{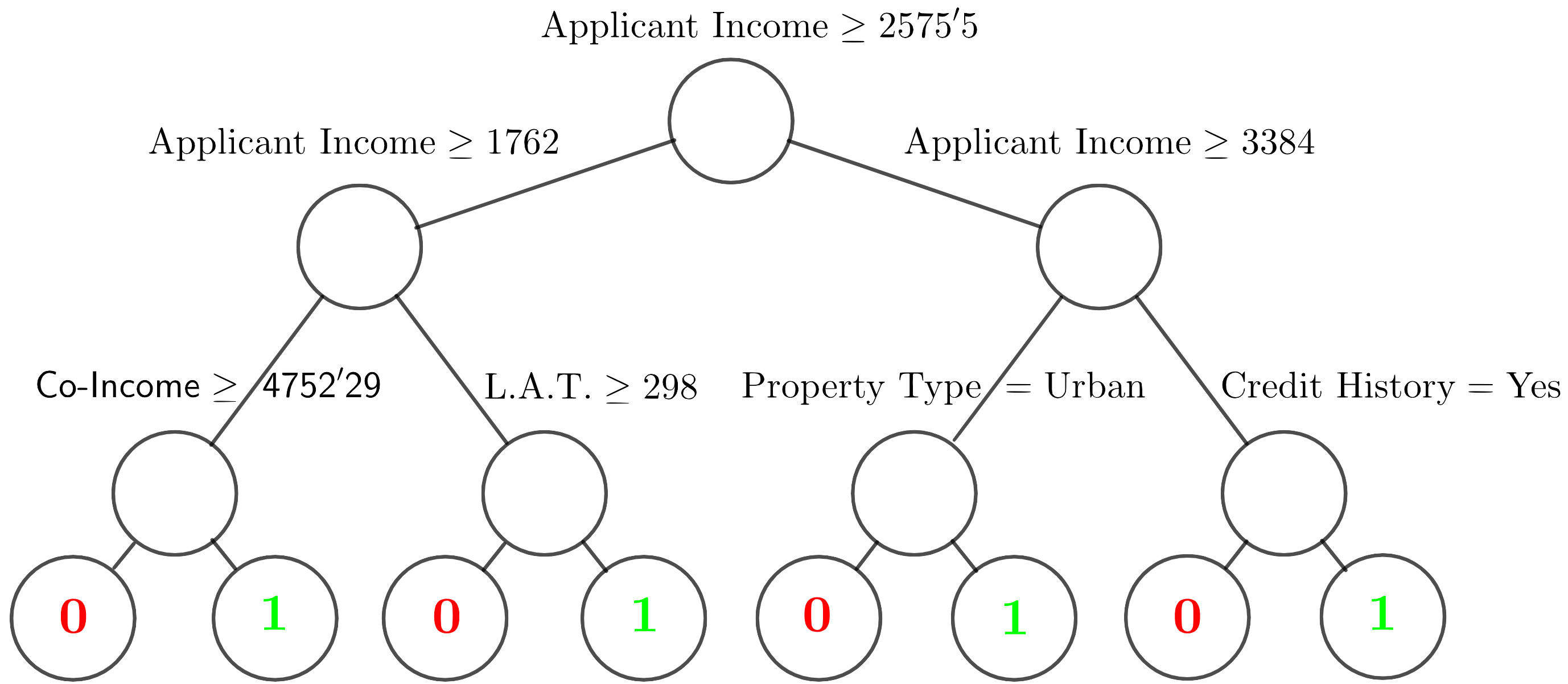}}
\caption{Solution OCT ($\sim 78\%$ of acccuracy) for the Mortgage case study.} \label{fig:4}
\end{figure}

\begin{figure}[h!]
\centering \fbox{\includegraphics[width=0.32\textwidth]{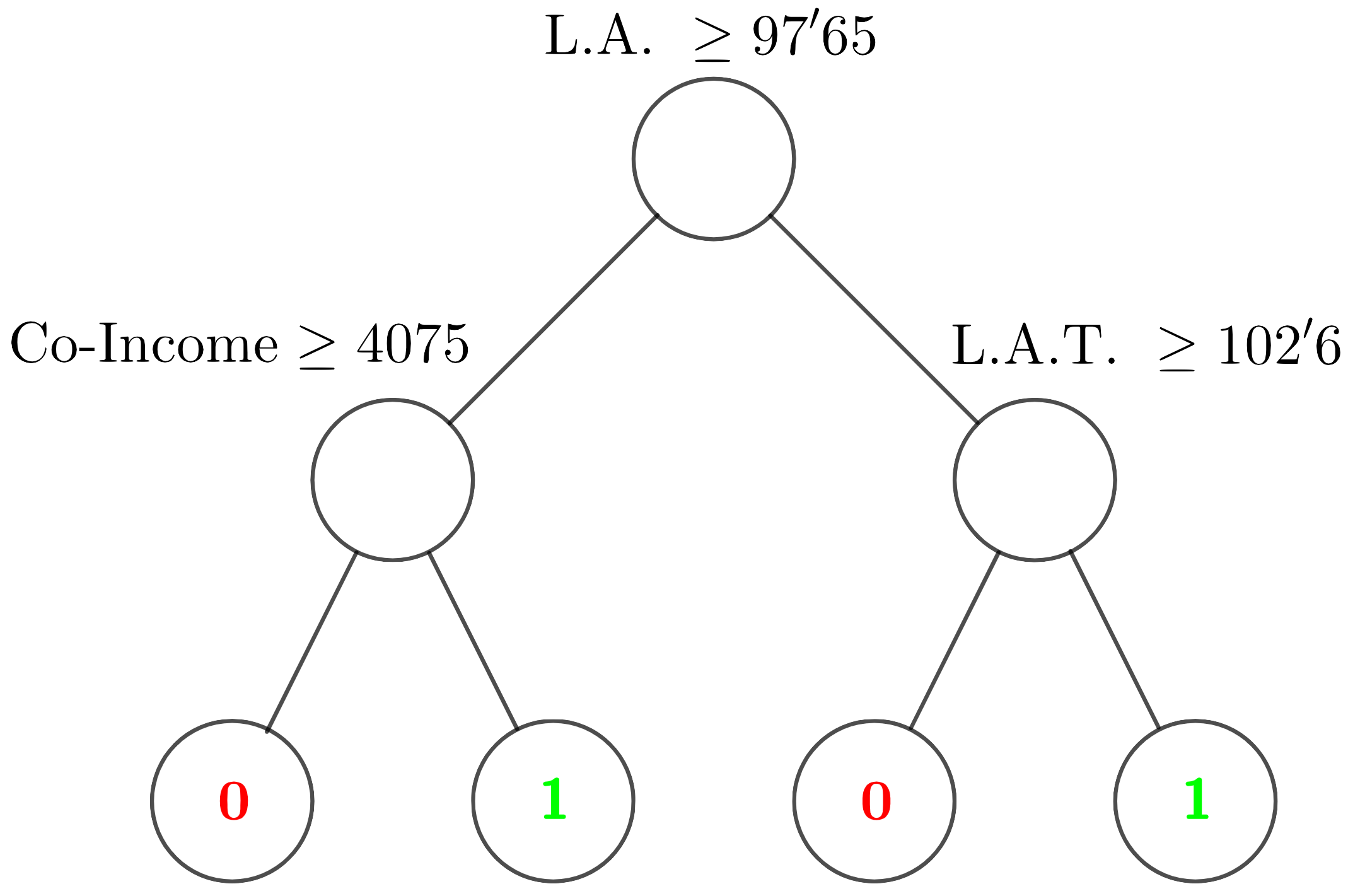}}~
\fbox{\includegraphics[width=0.274\textwidth]{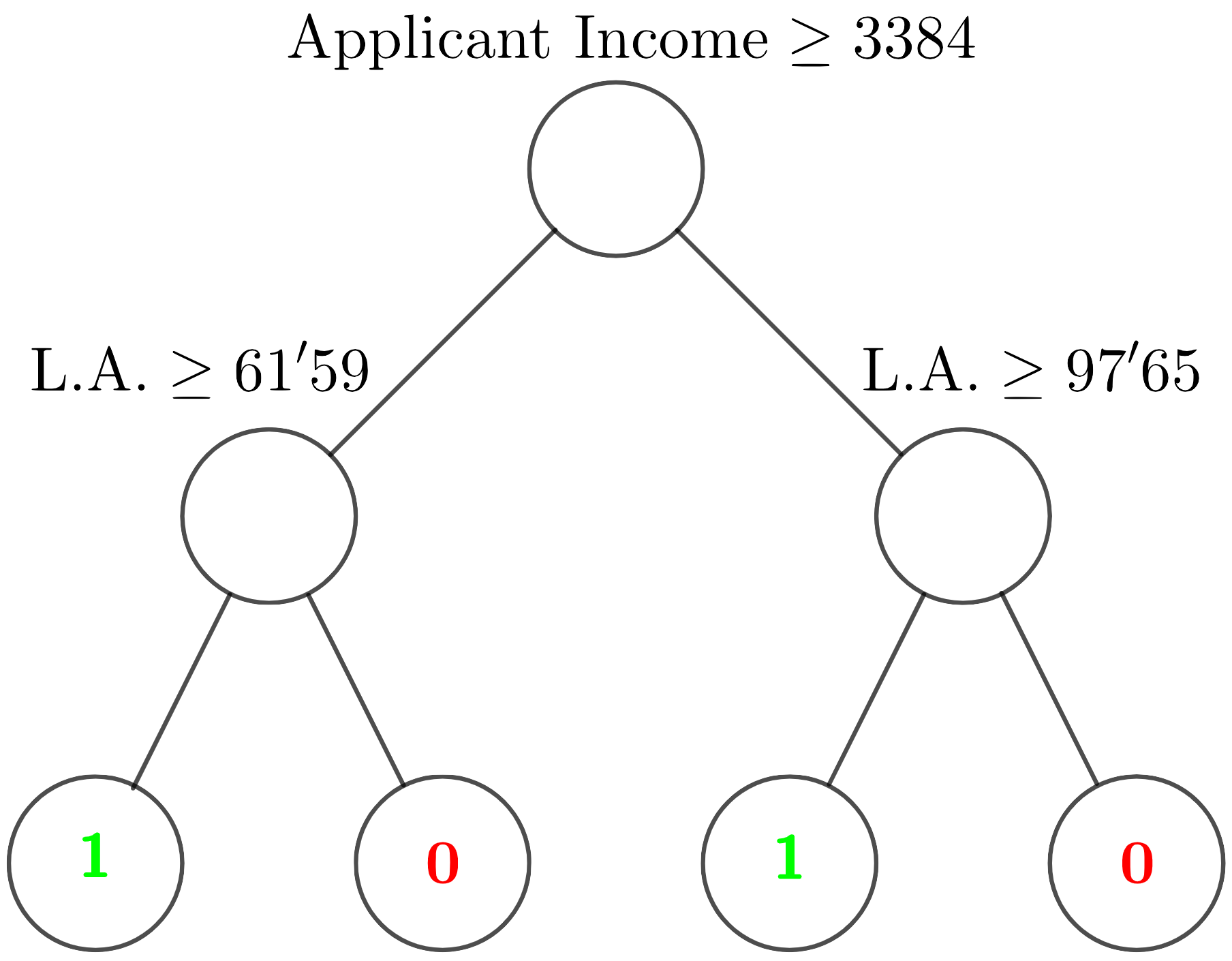}}~
\fbox{\includegraphics[width=0.32\textwidth]{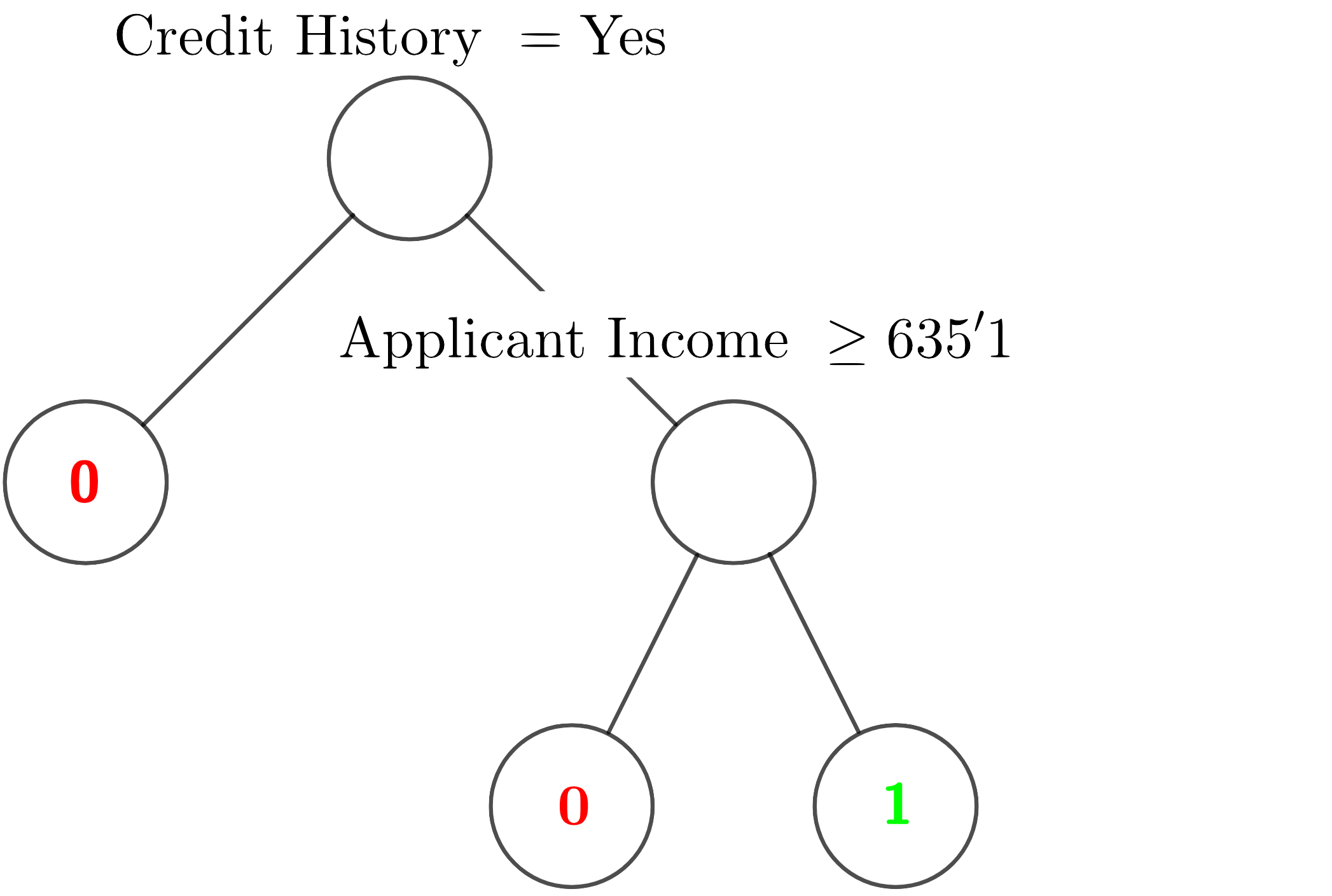}}
\caption{Solution 3-OCF ($\sim 82\%$ of accuracy) for the Mortgage case study.} \label{fig:5}
\end{figure}

In the following, we report, graphically, the solutions of an OCT of depth 3 and a 3-OCF of depth 2 to the problem. Given a partition of our dataset in train($50\%$)-validation($25\%$)-test($25\%$), we choose the model that works better in the validation set. The OCT is drawn in Figure \ref{fig:4}, whereas the WOCF is drawn in Figure \ref{fig:5}.

In the solution provided by OCT, the income variable is used in the first levels to derive the classification rule. In particular, the first split is built very close to the value of the first quartile of the income variable ($Q_1=2899$), and therefore, the L.A. and mainly L.A.T. variables, which is where we have the room for maneuver, will only affect to less than $25\%$ of the potential clients, while more than $75\%$ of potential clients will have no option to be recovered out from the ungranted set of clients. On the other hand, if we look at the solution provided by the 3-OCF, we find that one of the trees contains a branch that depends only on the variables L.A and L.A.T. Thus, the company will always have the option of calibrating these values to give a positive vote to a client, which, although will not ensure that the loan will be granted (the client must obtain at least $2$ out of $3$ votes to be granted), will increase its probability for all the clients. In fact, combining this branch with the credit history branch linked to having a minimum income allows for a lot of flexibility to end up granting a loan. In conclusion, in this case, the 3-OCF model not only provides a higher accuracy ($\sim 82\%$) than OCT ($\sim 78\%$), but also a more flexible solution for the company and the clients.

Although, as expected, the applicant income is a crucial variable to determine whether a client is granted, this is adequately combined with other variables (as L.A., L.A.T, Co-Income and Credit History) to alleviate negative values in one of the trees with positive values in the remainder trees. For instance, this might be the case of a client applying for a L.A. greater than 61.59, but with incomes smaller than $\$3385$. The second tree rejects the loan to this client. Nevertheless, in case the co-incomes are greater than $4075$, the first tree will accept the loan. Furthermore, if this client has a credit history and the client's incomes are greater than $635.1$, the third tree will also accept the loan. By majority, the client will get the loan. The solution obtained by OCT, for the same client with incomes smaller than $\$3384$, would grant the loan only if it is used to pay an urban property. Otherwise, the loan will not be accepted. 

\subsection{Case Study 2: Hotel Reservations}

In our second case study, we analyze one of the datasets provided by \cite{booking}. In such a work they present and study two datasets containing information on a series of room bookings in two hotels through the Booking company. In particular, we work with the first dataset which contains the information of a resort hotel, gathering $40060$ observations and $31$ variables, which cover all kinds of details about the customer's reservation. In {\color{blue}Table \ref{vars:hotel}}, we describe the variables that have been used in the models we have chosen  (the interested reader is referred to  \cite{booking} for further details on the variables defining the dataset).

\begin{table}[h!]
\centering
\begin{tabular}{lp{9.5cm}}
    Feature & Description\\\hline
    \texttt{Reservation Status} & Takes value 1 if the reservation is canceled and 0 otherwise.\\
    \texttt{Room Type} & Categorical variable.\\
    \texttt{ADR} & Average Daily Rate (total price paid by the customer divided by the number of nights).\\
    \texttt{Lead Time} & Number of days elapsed between the date of entry of the booking and the date of arrival.\\
    \texttt{TSR} & Total number of Special Requests of the room by the customer.\\
    \texttt{Parking} & Takes value 1 if the customer asks for a parking space and 0 otherwise.\\\hline
\end{tabular}
\caption{Predictive variables for the Hotel Reservations Dataset.\label{vars:hotel}}
\end{table}
The goal of this study is to decide whether a customer is susceptible to cancel a hotel reservation or not. In this way, the hotel desires to construct a  classification rule to determine the reliability of a customer when making the reservation. As a consequence of using a transparent classification rule, the hotels can detect the main reasons behind unstable bookings, and then derive the adequate interventions for more robust customers.

However, as in the previous case study, once a classification model has been trained, the hotel has some room for tweaking some of the variables, such as the room price, and could thus obtain more robust bookings from its customers (for example, a better adjustment on the pricing could cause a customer to no longer be classified as a potential cancellation but as a potential booking).\\

In our problem, the variable Room Type has not been used as a predictor variable but as a filter one. We made this decision because of its categorical nature and decided to focus on the $966$ type-6 room reservations. These rooms have on average the highest average cost per night, and also the highest imbalance in the target variable, i.e., the highest percentage of cancellations ($45\%$). Next, as done in the previous case study, we present the solutions provided by $3$-OCF and OCT for the problem, which have been obtained as in the previous case by following the train-validation-test partition schedule.

\begin{figure}[h!]
\centering
\fbox{\includegraphics[scale=0.7]{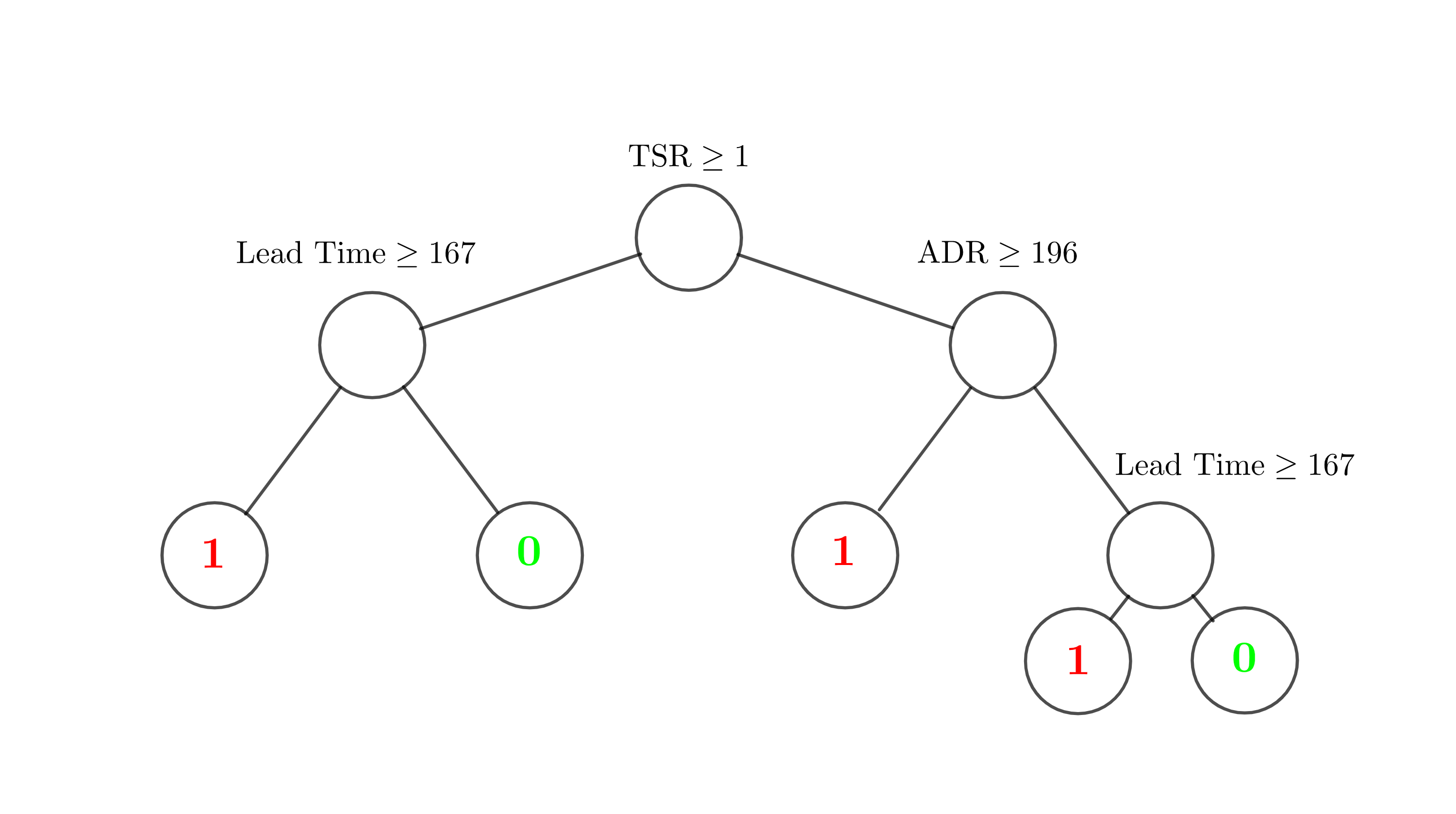}}
\caption{Solution OCT ($\sim 89\%$ of accuracy) for the Booking dataset.} \label{fig:6}
\end{figure}

\begin{figure}[h!]
\centering
\fbox{\includegraphics[scale=0.898]{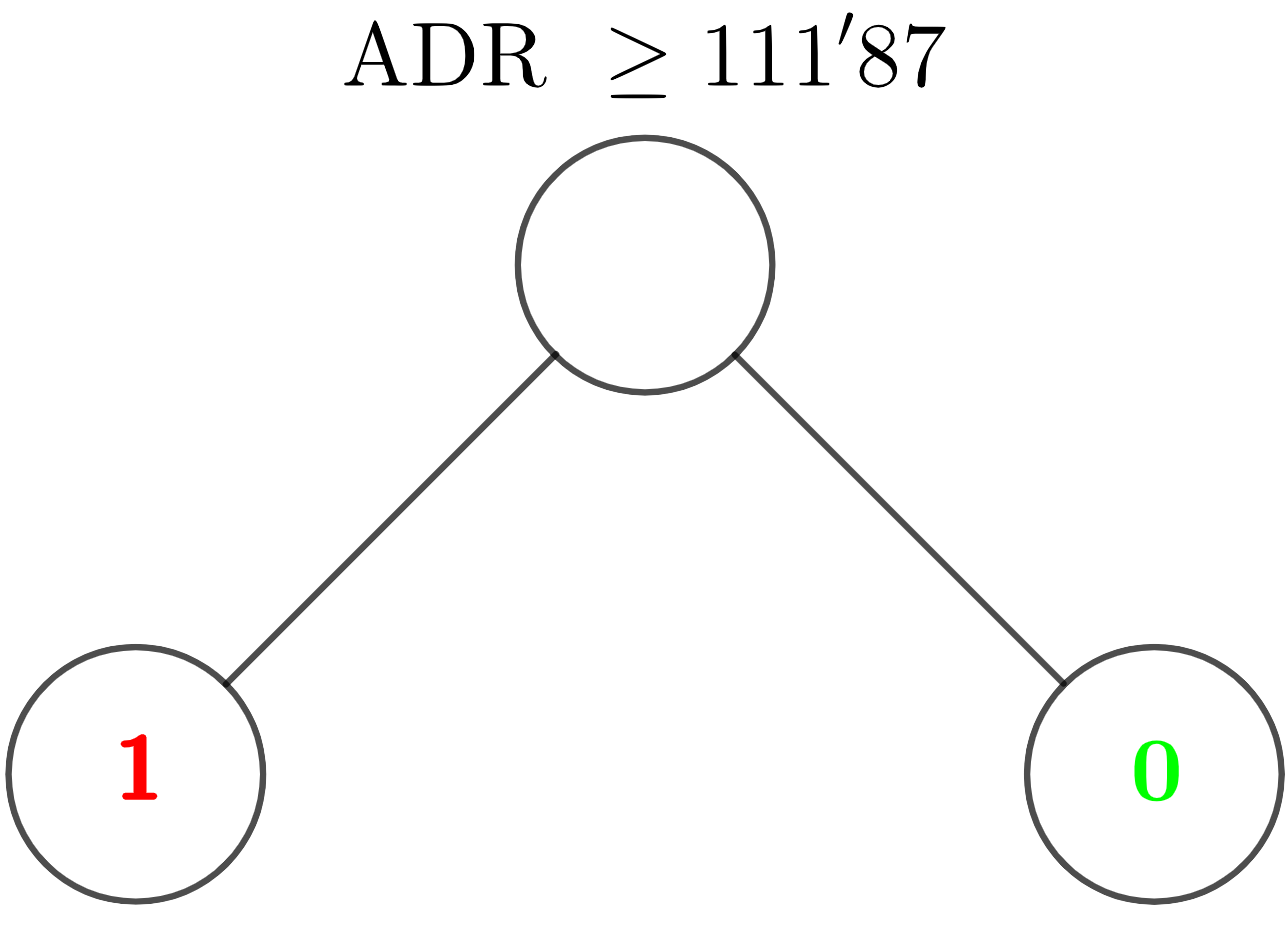}}~
\fbox{\includegraphics[scale=0.9]{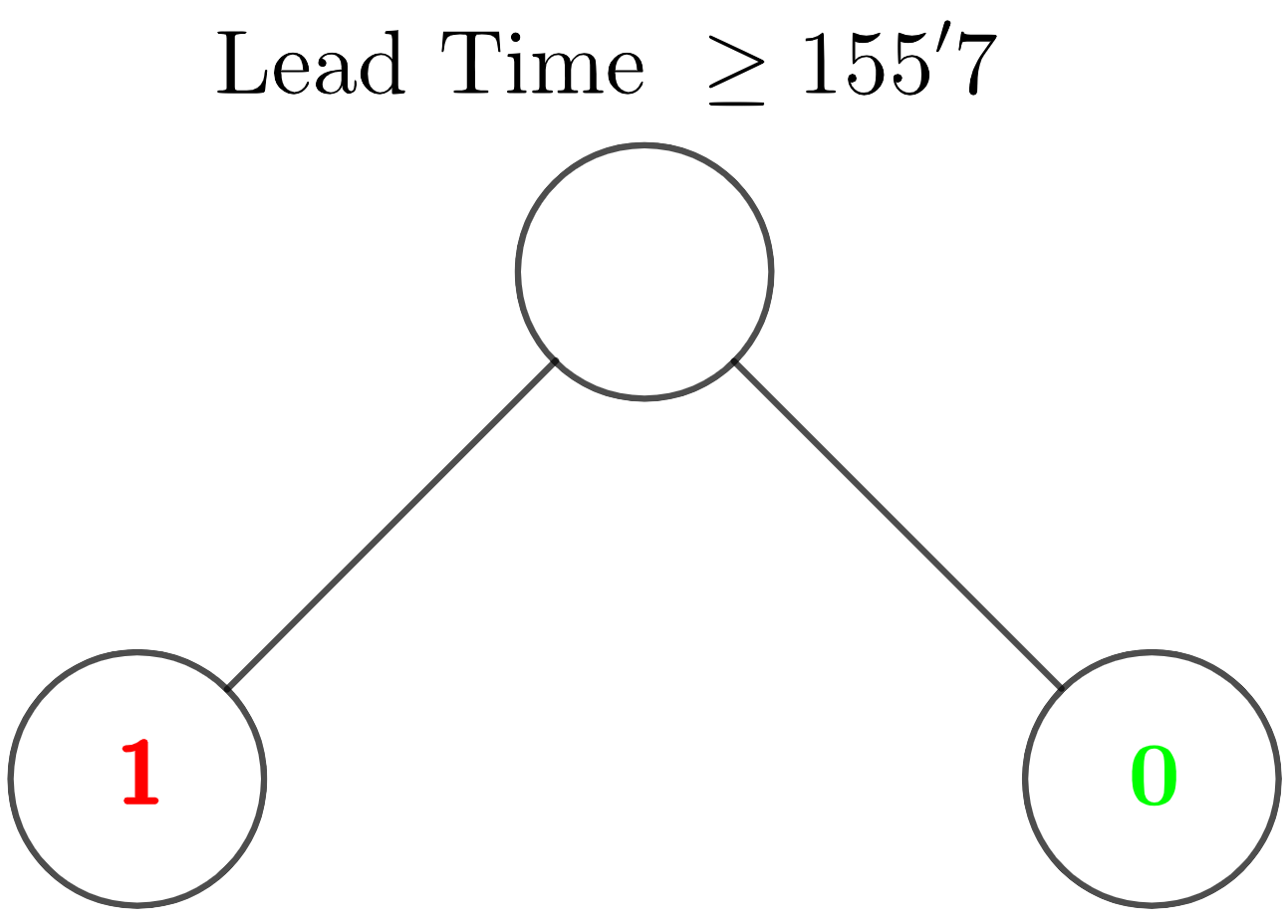}}~
\fbox{\includegraphics[scale=0.577]{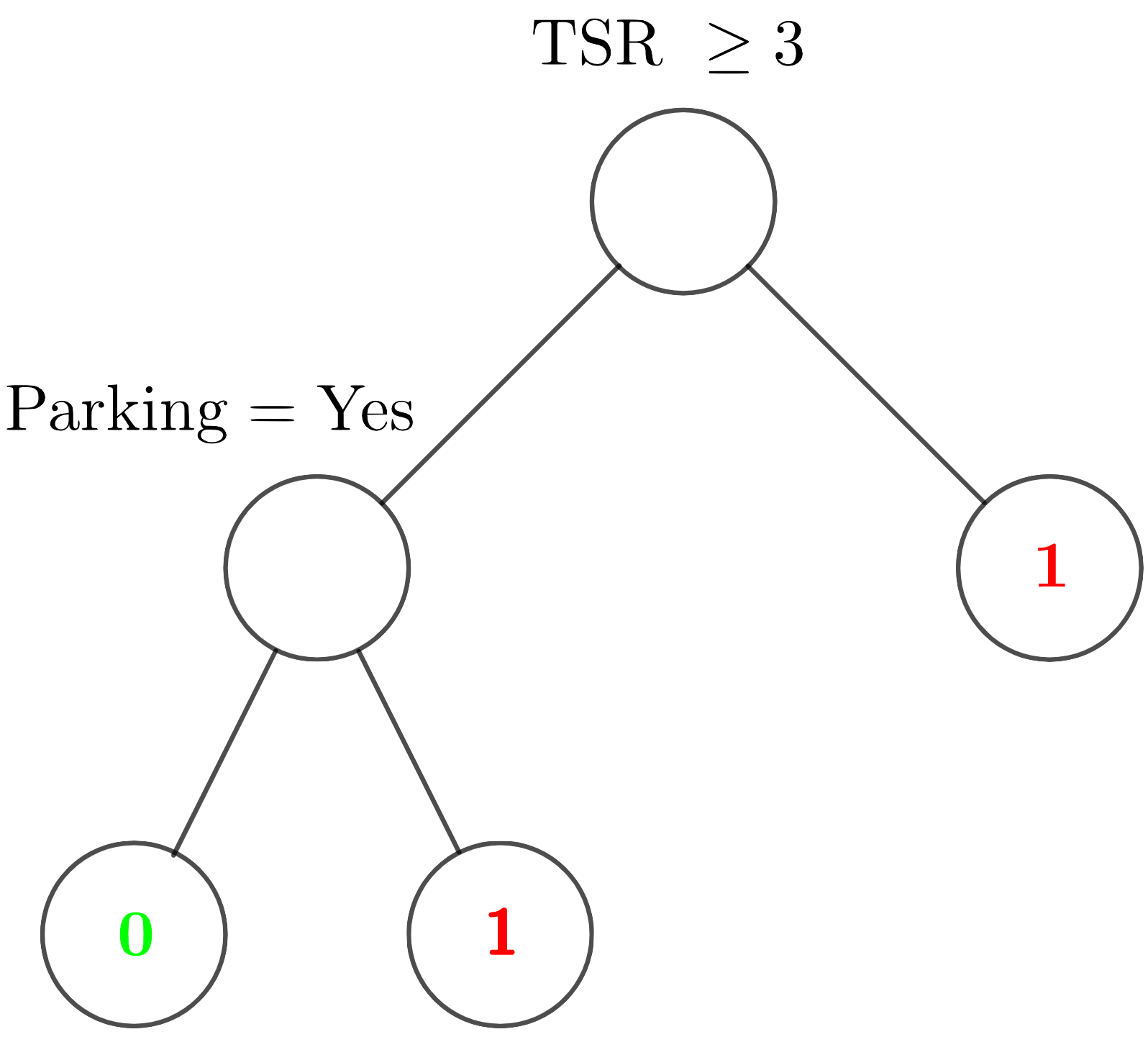}}
\caption{Solution 3-OCF ($\sim 89\%$ of accuracy) for the Bookng dataset.} \label{fig:7}
\end{figure}

In this case, both models have exactly the same performance in terms of accuracy in the test set ($\sim 89\%$) and both are easily interpretable (see Figures \ref{fig:6} and \ref{fig:7}). The OCT model starts by differentiating between the most demanding customers (given that we are dealing with the most expensive rooms in the hotel, it is to be expected that guests may have more specific requirements), and then takes into consideration the variables of price and lead time. In contrast, the 3-OCF model separates the price and lead time variables independently into two trees, and leaves a third tree to differentiate the most demanding customers, making a small distinction between those who ask for parking and those who do not (note that asking for parking implies that the customer has independence of movement and is, therefore, more likely to cancel and stay in other locations).\\

Thus, 3-OCF only automatically discards customers who book at a short advance and have a high number of special requests, whereas for the rest it will be able to adjust the price and maximize the hotel's bookings. On the other hand, the price variable in the OCT model can only be readjusted for a smaller subset of potential customers. Therefore, the $3$-OCF model provides more flexibility for the hotel to manage bookings successfully.

\subsection{Case Study 3: Airlines Satisfaction Problem}

In the previous two cases, we have seen the advantages that can be gained in terms of interpretability and flexibility by using 3-OCF for small/medium data sets. In this last case, we  show that this methodology can also be applied to large datasets and still obtains advantageous results in interpretability. In particular, we will analyze one of Kaggle's public datasets: \href{https://www.kaggle.com/datasets/teejmahal20/airline-passenger-satisfaction}{\texttt{airlines\_satisfaction}}.  The page itself invites to use of a set of $103,904$ observations and $25$ variables as training validation, and another $30,000$ observations as a test. We have used two-thirds of the training-validation observations. The remaining third of the observations have been used for validation while maintaining the proposed test observations.

The goal here is to determine, from the given dataset, whether a person is satisfied with an air travel experience ($44\%$ of the training sample) versus being dissatisfied or neutral ($56\%$ of the training sample). The $24$ predictor variables provide information about the customer and the flight itself. In Table \ref{vars:airlines} we describe the variables that have been involved in our chosen models, the reader can find the  details of the other variables in  \href{https://www.kaggle.com/datasets/teejmahal20/airline-passenger-satisfaction}{\texttt{airlines\_satisfaction}}.

\begin{table}[h!]
\centering
\begin{tabular}{lp{9.5cm}}
    Feature & Description\\\hline
 \texttt{Satisfaction} & takes value 1 if the customer is satisfied with the travel and 0 otherwise.\\
 \texttt{Online Boarding} & satisfaction with respect to the online boarding (from 0 to 5).\\
 \texttt{Wifi} & satisfaction with respect to the inflight wifi service (from 0 to 5).\\
 \texttt{Service} & satisfaction with respect to the inflight service (from 0 to 5).\\
 \texttt{Gate} & satisfaction with respect to the gate location (from 0 to 5).\\
 \texttt{Age} & Continuous variable.\\
 \texttt{Type} & Business / Personal.\\\hline
\end{tabular}
\caption{Predictive variables for the Hotel Reservations Dataset.\label{vars:airlines}}
\end{table}

\begin{figure}[h!]
\centering
\fbox{\includegraphics[scale=0.9]{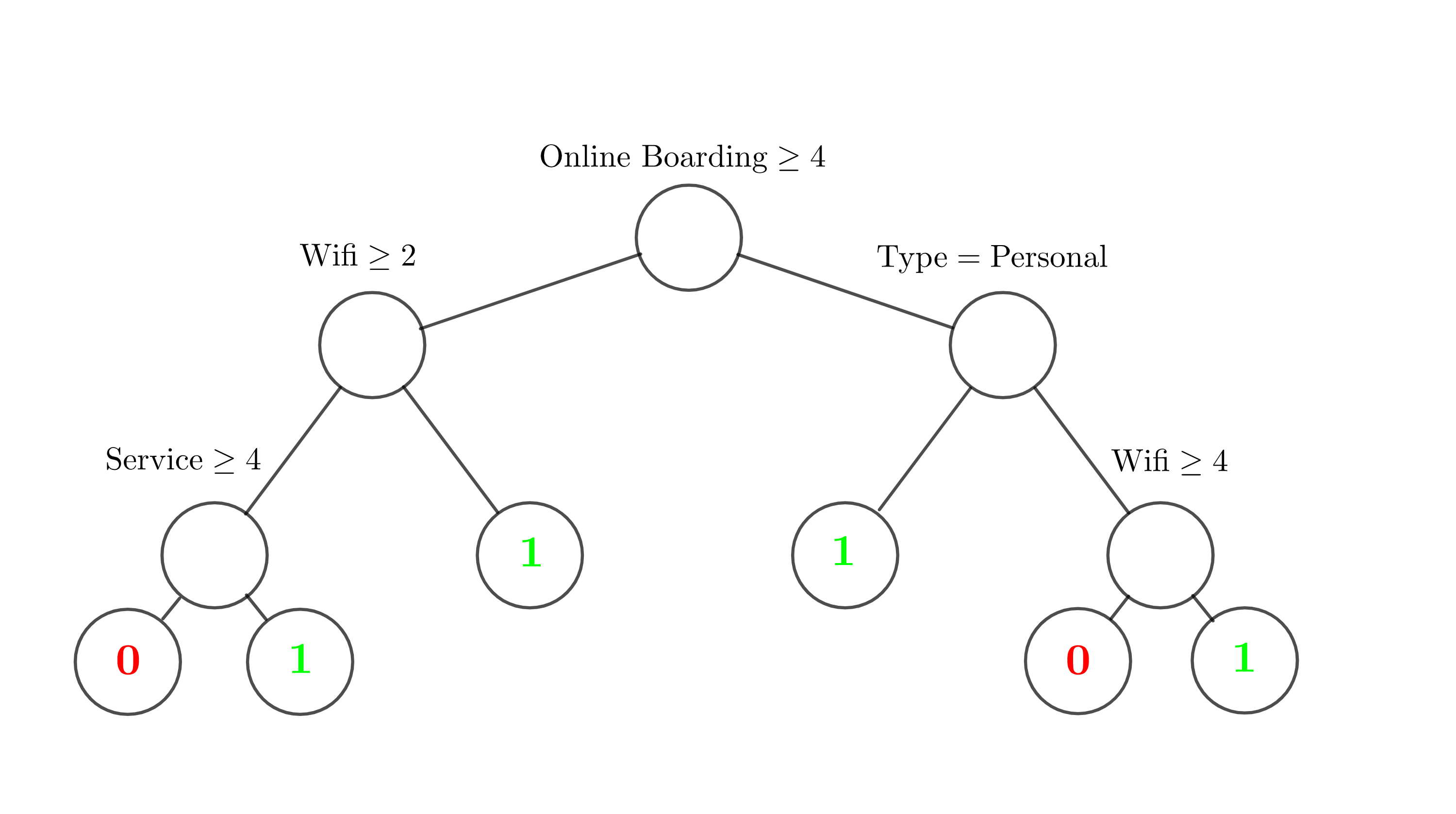}}
\caption{Solution OCT ($\sim 89\%$ of accuracy) for the airline dataset.} \label{fig:8}
\end{figure}

\begin{figure}[h!]
\centering
\fbox{\includegraphics[scale=0.64]{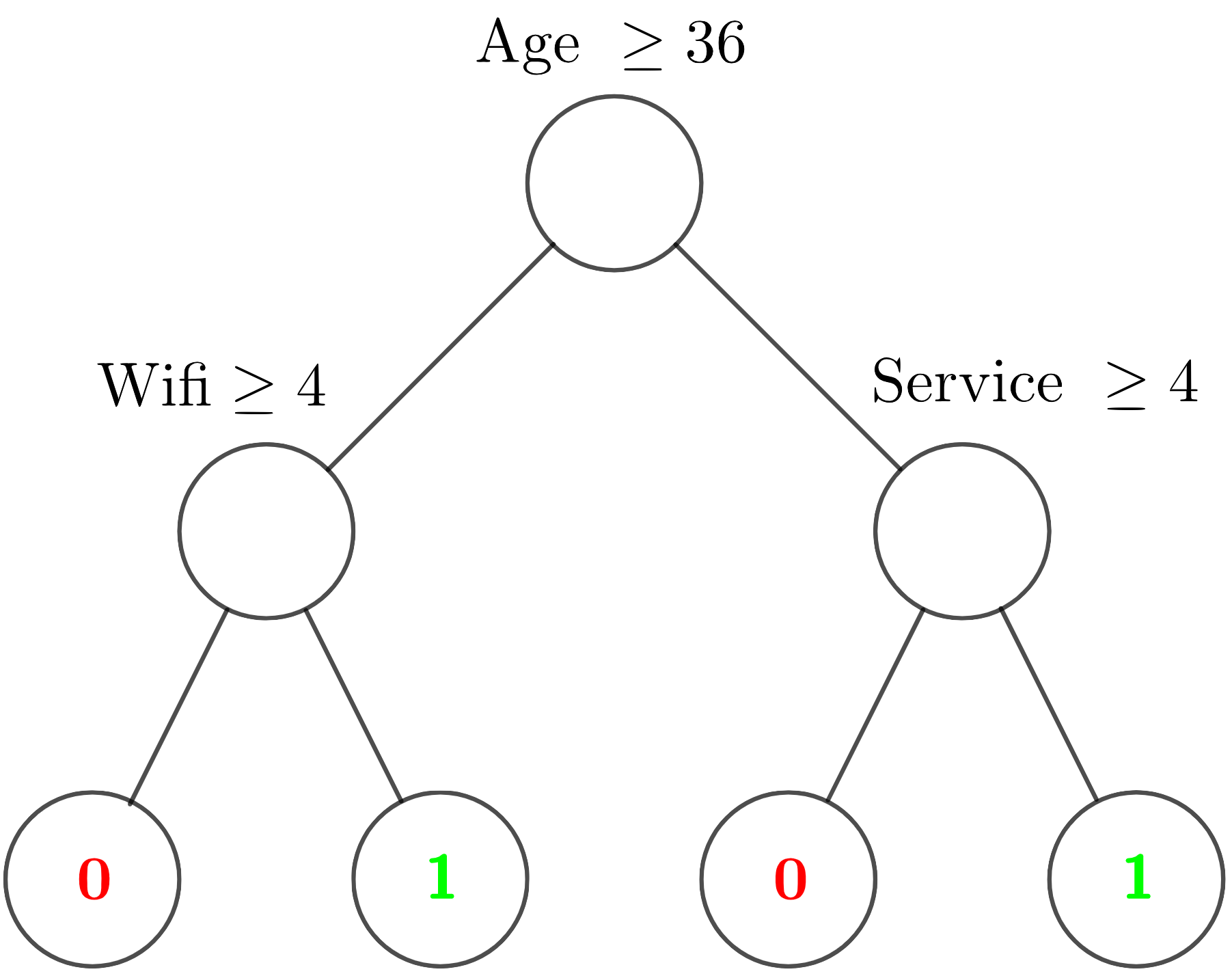}}~
\fbox{\includegraphics[scale=0.65]{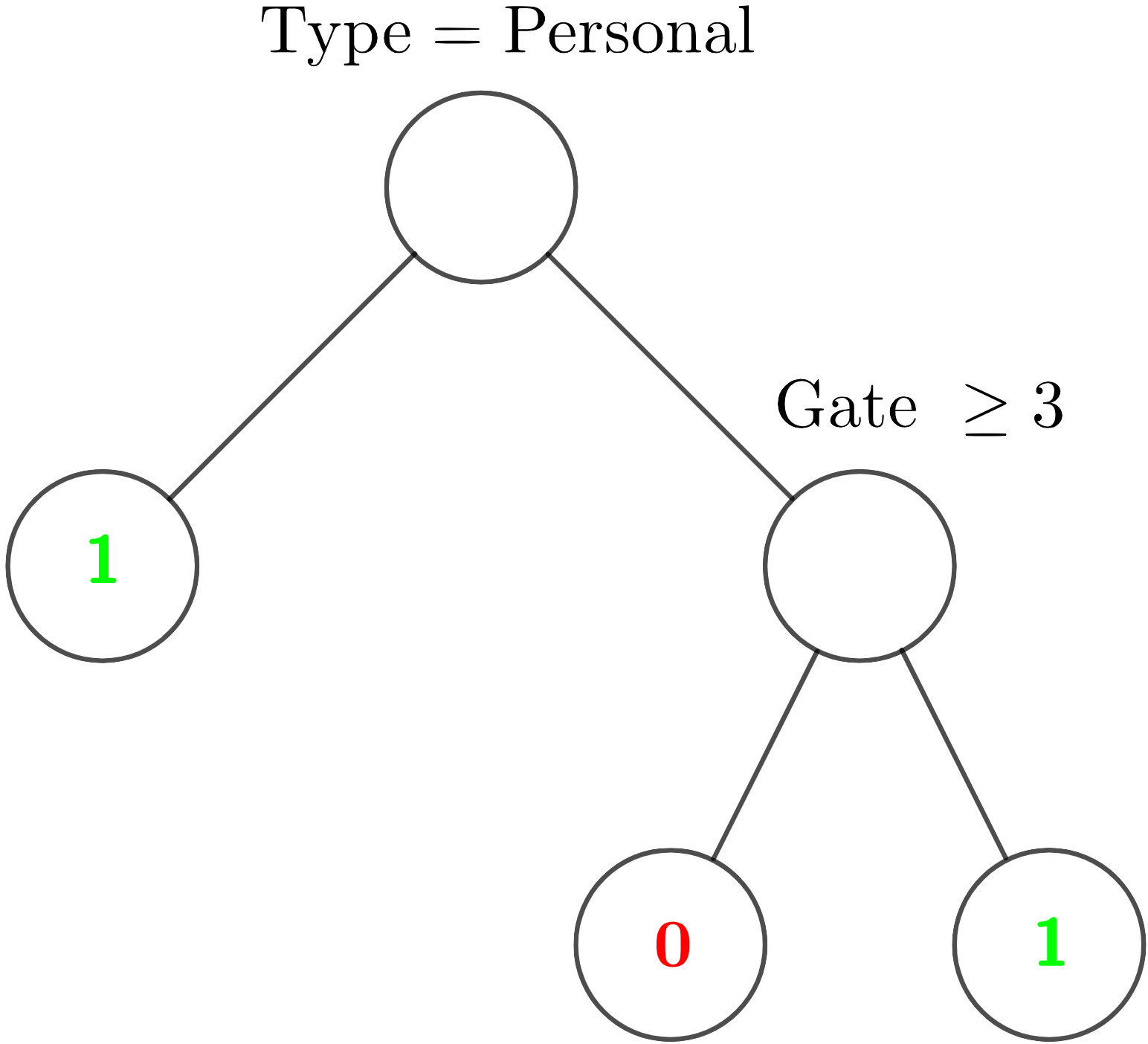}}~
\fbox{\includegraphics[scale=0.995]{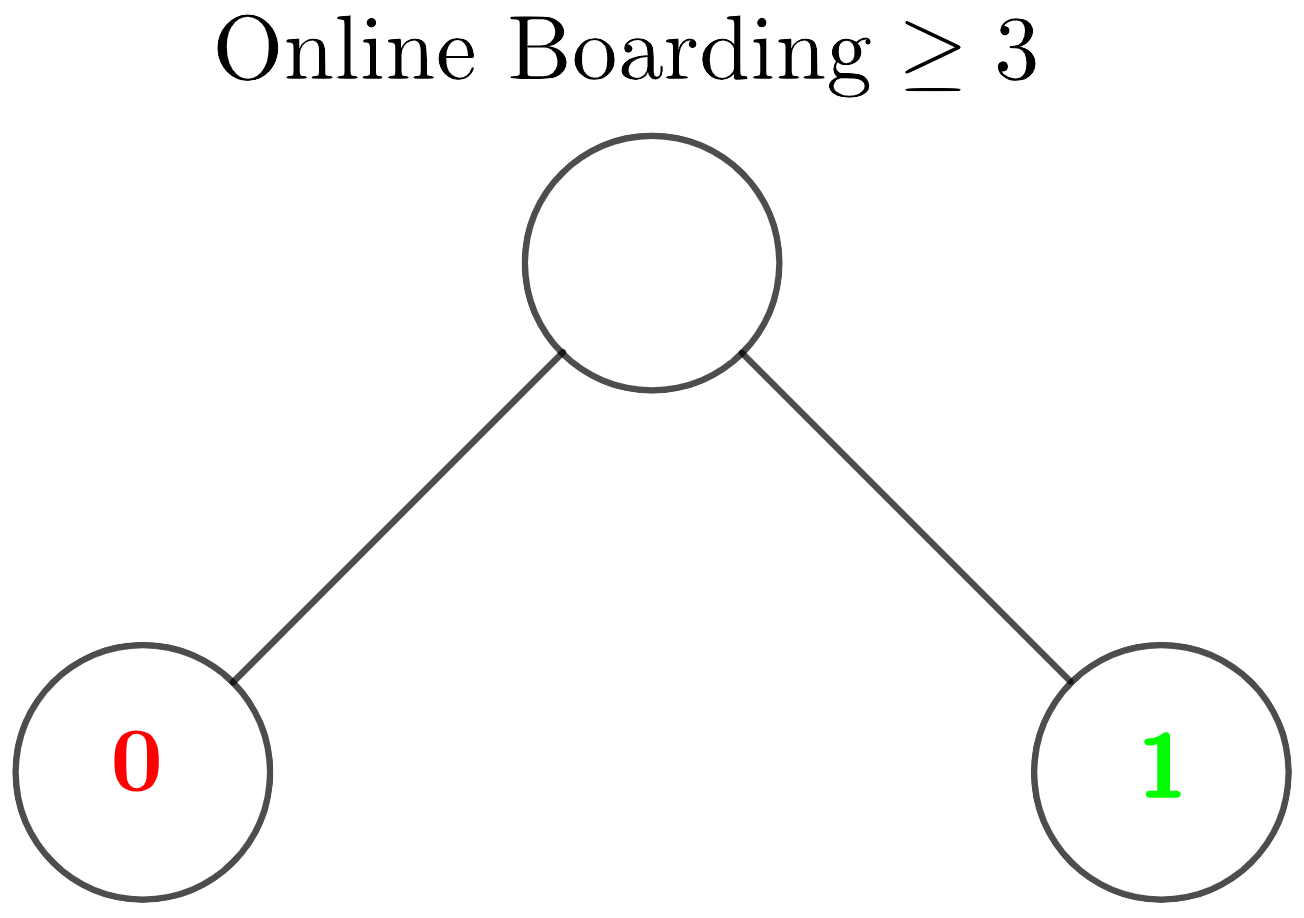}}
\caption{Solution 3-OCF ($\sim 89\%$ of accuracy) for the airline dataset.} \label{fig:9}
\end{figure}

The solutions obtained for this dataset are graphically shown in figures \ref{fig:8} (OCT) and \ref{fig:9} (WOCF). As can be observed, both methods obtain a similar accuracy in the test set (around $89\%$). Interestingly, in this case, the three trees obtained by 3-OCF provide information about three differentiated elements, that combined produce high accuracy. In the first tree, the first split differentiates customers by age. The younger customers rate positively  the trip if the plane has a good wifi system, whereas older customers, instead, take into account the inflight service to provide a positive rate. The second tree only gives a negative vote when the reason for the trip is personal and the customer is unsatisfied with the boarding gate, i.e., when the process of finding the gate and getting on the plane has not been smooth. In contrast, for business trips this is not a relevant factor and in general a positive vote is given. Finally, the third tree gives a positive vote if the online boarding system has been satisfactory. On the other hand, if we look at the solution proposed by OCT, although equally interpretable, the solution is not so clear at first sight. For example, we can see that the wifi variable is involved in different branches at different depths with different cut-off values, which may be less intuitive to analyze.
 
In conclusion, for this dataset, 3-OCF outputs three trees modeling different aspects of the satisfaction of airline trips, that seem more favorable than the one provided by the OCT, where following deeper branches with duplicated variables can make the interpretation a lot more confusing.

\section{Conclusions}\label{sec:conc} 

We propose here a new supervised classification method based on a voting-by-label process from a forest of classification trees. The decision rule is obtained by solving a mixed integer linear program on the training sample. The results have been tested on different benchmarking instances, outperforming some existing tree-based methods (CART and OCT) in terms of accuracy. With our approach we obtain results competitive to other tree ensemble based methods as Random Forests or XGB with a much larger number of trees by using a small amount of them. Finally, we show through three real case studies how WOCF solutions are not only easily interpretable but also more flexible than those provided by other classification tree methododologies.

Our proposed model is a first step in the study of this methodology, and a promising framework for future research. This is not a model that is intended to replace other methodologies that have proven to work well on classification problems, especially on large datasets. This is a new paradigm for researchers to add to their pool of tools when dealing with classification problems that require a high rate of accuracy as well as high explainability in moderate sized datasets. In these cases, it is often not easy to find an adequate tool to address the problem, and the use of one or another can have a major impact on the result. In view of the results, there is still room for improvement of WOCF in different directions. A deep analysis of the mathematical optimization problem that we propose would lead us to enhance its computational performance by finding valid inequalities, deriving strategies for fixing variables of the model, or designing decomposition approaches for solving larger-size training instances with less computational time. The study of alternative  mathematical optimization formulations for specially-shaped datasets, as those with binary or categorical features, will be also an interesting topic for a forthcoming paper.

\section*{Acknowledgements}

This research has been partially supported by Spanish Ministerio de Ciencia e Innovación, AEI/FEDER grant number PID 2020 - 114594GBC21 and Junta de Andalucía projects B-FQM-322-UGR20, C-EXP-139-UGR23, and AT 21\_00032. The first author was also partially supported by the IMAG-Maria de Maeztu grant CEX2020-001105-M /AEI /10.13039/501100011033 and UE-NextGenerationEU (ayudas de movilidad para la recualificaci\'on del profesorado universitario). This project is funded in part by Carnegie Mellon University’s Mobility21 National University Transportation Center, which is sponsored by the US Department of Transportation.

\end{document}